\documentclass[12pt]{amsart}
\usepackage[margin=0.9in]{geometry}
\usepackage{amsthm,amssymb,amsmath,xypic}
\usepackage[usenames,dvipsnames]{color}
\usepackage[colorlinks=true,citecolor=blue]{hyperref}
\usepackage{epigraph}
\usepackage{graphicx}
\usepackage{fancyhdr} 
\usepackage{enumitem}
\usepackage{mathrsfs}
\usepackage{enumitem}
\usepackage{rotating}
\usepackage{pdflscape}
\usepackage{caption}
\usepackage{subcaption} \newcommand{\be}{\begin{equation}}
\newcommand{\ee}{\end{equation}}
\newcommand{\bes}{\begin{equation*}}
\newcommand{\ees}{\end{equation*}}
\newcommand{\bea}{\begin{eqnarray}}
\newcommand{\eea}{\end{eqnarray}}
\newcommand{\beas}{\begin{eqnarray}}
\newcommand{\eeas}{\end{eqnarray}}
\newcommand{\ben}{\begin{note}}
\newcommand{\een}{\end{note}}
\newcommand{\bexl}{\vskip0.1em\noindent\hrulefill\vskip1em\begin{ExerciseList}}
\newcommand{\eexl}{\end{ExerciseList}\hrulefill}

\newcommand{\bthm}{\begin{theorem}}
\newcommand{\ethm}{\end{theorem}}
\newcommand{\bpro}{\begin{prop}}
\newcommand{\epro}{\end{prop}}
\newcommand{\bcor}{\begin{corollary}}
\newcommand{\ecor}{\end{corollary}}
\newcommand{\bcon}{\begin{conjecture}}
\newcommand{\econ}{\end{conjecture}}
\newcommand{\bp}{\begin{proof}}
\newcommand{\ep}{\end{proof}}
\newcommand{\blem}{\begin{lemma}}
\newcommand{\elem}{\end{lemma}}
\newcommand{\bn}{\begin{note}}
\newcommand{\en}{\end{note}}
\newcommand{\benum}{\begin{enumerate}}
\newcommand{\eenum}{\end{enumerate}}
\newcommand{\bed}{\begin{defn}}
\newcommand{\eed}{\end{defn}}
\newcommand{\brem}{\begin{remark}}
\newcommand{\erem}{\end{remark}}

\newcommand{\btik}{\begin{tikzpicture}\begin{axis}[scale=0.5,axis y line=center, axis x line=middle]}
\newcommand{\etik}{\end{axis}\end{tikzpicture}}

\let\mapsto=\longmapsto

\newcommand{\upperRomannumeral}[1]{\uppercase\expandafter{\romannumeral#1}}

 \usepackage{ifthen}
\usepackage{natbib}
\let\cite=\citep
\usepackage{zref-clever}
\let\Cref=\zcref
\usepackage{afterpage}
\newtheorem{theorem}[equation]{Theorem}      \newtheorem{lemma}[equation]{Lemma}          \newtheorem{corollary}[equation]{Corollary}  \newtheorem{proposition}[equation]{Proposition}

\theoremstyle{definition}
\newtheorem{conj}[equation]{Conjecture}
\newtheorem{example}[equation]{Example}

\theoremstyle{definition}
\newtheorem{defn}[equation]{Definition}
\theoremstyle{remark}

\theoremstyle{definition}
\newtheorem{remark}[equation]{Remark}

\numberwithin{equation}{section}

\newcommand{\fixnumberwithin}[1]{
\numberwithin{equation}{#1}
	\numberwithin{theorem}{#1}
	\numberwithin{conj}{#1}
	\numberwithin{lemma}{#1}
	\numberwithin{proposition}{#1}
	\numberwithin{corollary}{#1}
	\numberwithin{defn}{#1}
	\numberwithin{remark}{#1}
	\numberwithin{rem}{#1}
	\numberwithin{question}{#1}
\IfPackageLoadedTF{zref-clever}{
		\ifundef{\Cref}{\let\Cref=\zcref}{\relax} \ifundef{\cref}{\let\cref=\zcref}{\relax} \zcsetup{cap} \zcRefTypeSetup{conj}{
			Name-sg = Conjecture ,
			name-sg = conjecture ,
			Name-pl = Conjectures ,
			name-pl = conjectures ,
		}
		\zcRefTypeSetup{defn}{
			Name-sg = Definition,
			name-sg = definition ,
			Name-pl = Definitions ,
			name-pl = definitions ,
		}
	}
	{\relax}
}

\zcsetup{cap} \AddToHook{env/theorem/begin}{\zcsetup{countertype={equation=theorem}}}
\AddToHook{env/lemma/begin}{\zcsetup{countertype={equation=lemma}}}
\AddToHook{env/proposition/begin}{\zcsetup{countertype={equation=proposition}}}
\AddToHook{env/theoremdef/begin}{\zcsetup{countertype={equation=theoremdef}}}
\AddToHook{env/defn/begin}{\zcsetup{countertype={equation=defn}}}
\AddToHook{env/remark/begin}{\zcsetup{countertype={equation=remark}}}
\AddToHook{env/rem/begin}{\zcsetup{countertype={equation=rem}}}
\AddToHook{env/conj/begin}{\zcsetup{countertype={equation=conj}}}
\AddToHook{env/question/begin}{\zcsetup{countertype={equation=question}}}
\AddToHook{env/example/begin}{\zcsetup{countertype={equation=example}}}

\zcRefTypeSetup{equation}{
	Name-sg = eq. ,
	name-sg = eq. ,
	Name-pl = eqns. ,
	name-pl = eqns. ,
}

\zcRefTypeSetup{theorem}{
	Name-sg = Theorem ,
	name-sg = theorem ,
	Name-pl = Theorems ,
	name-pl = theorems ,
}

\zcRefTypeSetup{defn}{
	Name-sg = Definition ,
	name-sg = definition ,
	Name-pl = Definitions ,
	name-pl = definitions ,
}

\zcRefTypeSetup{theoremdef}{
	Name-sg = Theorem-Definition ,
	name-sg = theorem-definition ,
	Name-pl = Theorem-Definitions ,
	name-pl = theorem-definitions ,
}

\zcRefTypeSetup{proposition}{
	Name-sg = Proposition ,
	name-sg = proposition ,
	Name-pl = Propositions ,
	name-pl = propositions ,
}

\zcRefTypeSetup{lemma}{
	Name-sg = Lemma ,
	name-sg = lemma ,
	Name-pl = Lemmas ,
	name-pl = lemmas ,
}

\zcRefTypeSetup{remark}{
	Name-sg = Remark ,
	name-sg = remark ,
	Name-pl = Remarks ,
	name-pl = remarks ,
}

\zcRefTypeSetup{question}{
	Name-sg = Question ,
	name-sg = question ,
	Name-pl = Questions ,
	name-pl = questions ,
}

\zcRefTypeSetup{example}{
	Name-sg = Example ,
	name-sg = example ,
	Name-pl = Examples ,
	name-pl = examples ,
}

\zcRefTypeSetup{conj}{
	Name-sg = Conjecture ,
	name-sg = conjecture ,
	Name-pl = Conjectures ,
	name-pl = conjectures ,
}

\let\cite=\citep 

\let\isom=\simeq

\newcommand{\abs}[1]{\left\vert#1\right\vert}

\newcommand{\gal}{{\rm Gal}}

\newcommand{\Q}{{\mathbb Q}}

\newcommand{\Z}{{\mathbb Z}}

\renewcommand{\O}{{\mathcal O}}

\renewcommand{\bpro}{\begin{proposition}}
\renewcommand{\epro}{\end{proposition}}

\newcommand{\np}[1]{N_{#1}(E)}
\let\vare=\varepsilon

\newcounter{constantCounter}
\newcommand{\cnst}{c_{\arabic{constantCounter}\stepcounter{constantCounter}}}
\newcommand{\cnste}{c_{\arabic{constantCounter}}(\vare)\stepcounter{constantCounter}}
\newcommand{\labelCounter}[1]{\expandafter\def\csname#1\endcsname{c_{\arabic{constantCounter}}}}

\newcommand{\lx}{\log x}
\newcommand{\llx}{\log\log x}

\newcommand{\lix}{{\rm Li}(x)}
\newcommand{\npp}{\np{p}}
\newcommand{\setx}[1]{\mathcal{#1}(x)}
\newcommand{\csetx}[1]{|\setx{#1}|}
\newcommand{\ax}{\setx{A}}
\newcommand{\cax}{\csetx{A}}
\newcommand{\bx}[1]
{\ifthenelse
	{\equal{#1}{}}{\setx{B}}{{\mathcal B}_{#1}(x)}}
\newcommand{\cbx}[1]{|\bx{#1}|}
\newcommand{\bxx}[2]{{\mathcal B}_{#1}(#2)}
\newcommand{\cbxx}[2]{|\bxx{#1}{#2}|}
\newcommand{\cx}{\setx{C}}
\newcommand{\ccx}{\csetx{C}}
\newcommand{\dx}{\setx{D}}
\newcommand{\cdx}{\csetx{D}}

\newcommand{\low}{(1-\vare)\llx}

\newcommand{\dxo}{{\mathcal{D}}_1(x)}
\newcommand{\dxt}{{\mathcal{D}}_2(x)}
\newcommand{\cdxo}{|\dxo|}
\newcommand{\cdxt}{|\dxt|}

\newcommand{\sumtwo}{\sum\nolimits_2}

\newcommand{\sump}[1]{\sum_{p\leq{#1}}}
\newcommand{\grh}{Let $E/\Q$ be an elliptic curve. If $E$ does not have complex multiplication, then assume that GRH (\Cref{def:grh}) holds for $E$.}
\begin{document}

\title[]{On the equation $\np{p_1}\np{p_2}\cdots\np{p_k}=n$}\author{Kirti Joshi}\address{Math. department, University of Arizona, 617 N Santa Rita, Tucson
85721-0089, USA.}

\thanks{}\subjclass{}

\newcommand{\erdos}{Erd\H os}

\begin{abstract}
	For a given elliptic curve $E/\Q$, let $N_p(E)$  be the number of points on $E$ modulo $p$ for a prime of good reduction for $E$.  Given integer $n$, let $G_k(E,n)$ be the number of $k$-tuples of $p_1<p_2<\ldots <p_k$ primes of good reduction for $E$, for which the equation in the title holds, then on assuming the Generalized Riemann Hypothesis for elliptic curves without CM (and unconditionally if the curves have complex multiplication), I show that $\varlimsup_{n\to\infty} G_k(E,n)=\infty$ for any integer $k\geq 3$.  I conjecture that this result also holds for $k=1,2$ i.e.  this conjecture says that there are arbitrarily long ``elliptic progressions of primes'' i.e. sequences of primes $p_1<p_2<\cdots <p_m$ of arbitrary lengths $m$ such that $\np{p_1}=\np{p_2}=\cdots =\np{p_m}$.
\end{abstract}
\maketitle
\epigraph{N was a net\\
	Which was thrown in the sea\\
	To catch fish for dinner\\
	For you and for me.}{\citeauthor{lear}}
\tableofcontents

\section{Introduction}\label{se:intro}
Let me begin with an example which explains the problem I consider in this paper. 
Consider the elliptic curve $E:y^2+y=x^3-x$ over $\Q$. This is an elliptic curve with conductor $37$. Let $p$ be a prime of good reduction, and let $N_p(E)$ be the number of points on this elliptic curve modulo $p$.
It is a well-known theorem of Hasse \citep[Chapter V, Theorem 1.1]{silverman-arithmetic} that for all the  primes of good reduction $p$ one has
\be\label{eq:hw-bound}
p+1-2\sqrt{p}\leq \np{p} \leq p+1+2\sqrt{p}.
\ee

A simple computation with \citep{sage} and \citep{pari} reveals that following equalities hold:
\be\label{eq:elliptic-examples1}
\begin{aligned}\np{2}\cdot\np{13}\cdot\np{43}&= &3360\ \ \ \ \ \ \ &=
\np{3}\cdot\np{5}\cdot\np{67}\\
\np{5}\cdot\np{43}\cdot\np{73}&=&25200\ \ \ \ \ \ &=\np{17}\cdot\np{19}\cdot\np{61}\\
\np{101}\cdot\np{107}\cdot\np{251}&=&3107520\ \ \ \ \ &=\np{113}\cdot\np{127}\cdot\np{167}\\
\np{1009}\cdot\np{1181}\cdot\np{1601}&=&1988217000\ \ \ \ &=\np{1063}\cdot\np{1283}\cdot\np{1399}.
\end{aligned}
\ee
All  the triples of primes entering these equations are all pairwise distinct; moreover the six numbers which enter these equalities are typically distinct, but occasionally not. 

More generally, let $E/\Q$ be any elliptic curve and fix an integer $k\geq 1$. Given a natural number $n\geq 1$, I consider the equation:
\be\label{eq:erdos} 
\np{p_1}\cdots \np{p_k} = n.
\ee 
in pairwise distinct,  primes of good reduction $p_1,\ldots,p_k$.

For $k\geq 1$, let $G_k(E,n)$  be the number of solutions to \eqref{eq:erdos} in pairwise distinct  primes of good reduction $p_1<p_2<\cdots<p_k$ for $E$. 

The number $G_k(E,n)$ is always finite. To see this, I claim the lower bound 
\be\label{eq:triv-bound} N_p(E)\geq p/100 \qquad \text{for all primes } p.\ee 
Indeed,  from Hasse's lower bound \eqref{eq:hw-bound}, for any prime $p$, one has
\be 
\frac{N_p(E)}{p}\geq \frac{p+1-2\sqrt{p}}{p}\geq\begin{cases}
	1-\frac{2}{\sqrt{p}}\geq \frac{1}{100} & {\rm if}\ p\geq 5\\
	1+\frac{1}{3}-\frac{2}{\sqrt{3}}\geq \frac{1}{100} & {\rm if}\ p=3\\
	1+\frac{1}{2}-\frac{2}{\sqrt{2}}\geq \frac{1}{100} & {\rm if}\ p=2\\
\end{cases}
\ee
Hence  from \eqref{eq:erdos} and \eqref{eq:triv-bound} one sees that
\be n\geq  \frac{p_1\cdots p_k}{100^k}.\ee 
Thus, for any fixed integer $k\geq1$, the $k$-tuples of primes which contribute to $G_k(E,n)$ are finite in number. Hence for any elliptic curve $E/\Q$, any fixed integer $k\geq 1$, the number $G_k(E,n)$ is finite for any integer $n>1$. 

Throughout the paper, the phrase \textit{$E/\Q$ has complex multiplication (CM)}, will mean that $E/\Q$ is an elliptic curve with complex multiplication by the ring of integers $\O_K$ of an imaginary quadratic field $K$.

The purpose of this paper is to prove  the following theorem which shows that  these equalities \eqref{eq:elliptic-examples1} are not small numerical accidents. 
\bthm\label{th:main}
\grh\ 
Then for any integer $k\geq 3$
$$ 
	\varlimsup_{n\to\infty} G_k(E,n)=\infty.
$$
More precisely, for every integer $k\geq 3$ and for all sufficiently large $x$, there exist  integers $n\leq x$ with $$G_k(E,n)\geq (\lx)^{\delta}$$
and $\delta=\delta(k)= \log(2)k(1-2\varepsilon)-2>0$.
\ethm

In the CM case, I show in \Cref{cor:cm-case}, that a stronger assertion holds, namely \Cref{th:main} is true for $k\geq 3$ even if the primes in \eqref{eq:erdos} are required to be primes of good ordinary reduction.

For the definition of the term \textit{Generalized Riemann Hypothesis (GRH) holds  for $E/\Q$} see \Cref{def:grh}. Let me also point out that one can also replace the assumption that GRH holds for $E/\Q$ by other types of hypothesis such as  $\theta$-quasi-GRH etc. (\cite{cojocaru05}), however since I make no claim of providing sharpest error terms nor the best possible constants, the choice of stronger or weaker hypothesis is insignificant (to me). At any rate, it is unlikely that the growth of $G_k(E,n)$ given by the proof of \Cref{th:main} is optimal.

Let me point out that  in \citep{erdos36}, Paul \erdos\  considered the question of the number of solutions (in triples of pairwise distinct primes $p_1,p_2,p_3$) to the equation 
\be 
(p_1-1)(p_2-1)(p_3-1)=n,\ee 
and proved the corresponding result in this case 
and that his paper served as an inspiration for computations of the examples \eqref{eq:elliptic-examples1} which led to this note.  I am not aware of any paper of \erdos\ where he addresses the problem of the number of solutions in the general case $$(p_1-1)(p_2-1)\cdots (p_k-1)=n$$ for $k\geq2$, and hence this is proved here in \Cref{th:main2} using methods developed for the proof of \Cref{th:main}. In this classical case, the result is unconditional and is established for all $k\geq 2$ ($k=2$ falls back on \cite{erdos36b} which treats the problem of counting solutions  $n=(p-1)(q-1)$ in distinct primes $p,q$).

While \citep{erdos36} served as an inspiration for this paper,  my approach  differs from \erdos' approach in several important points. Notably, I use the fourth moment to prove \Cref{le:convergence} instead of the Hardy-Ramanujan Theorem. I also prove \Cref{th:erdos-product}, which was proved by \erdos\ for $k=3$ in \citep{erdos36}, a bit differently.   

More important observation of this paper is this: the case $k=1$ seems to be interesting in its own right: I conjecture (in \Cref{con:k-one-case}{\bf(1)}) that a result similar to \Cref{th:main} holds for $k=1$. Specifically, \Cref{con:k-one-case}{\bf(1)} asserts that there exists a sequence of integers $n\to\infty$ and a sequence of primes $p_1<p_2<\cdots<p_{m_n}$ with $m_n\to\infty$ (with $n$) such that 
\be\label{eq:elliptic-prog}\np{p_1}=\np{p_2}=\cdots=\np{p_{m_n}}=n.\ee 
I call the tuple  $p_1<p_2<\cdots<p_{m_n}$ of primes satisfying \eqref{eq:elliptic-prog} an \textit{elliptic progression of primes (of length $m_n$) given by the curve $E/\Q$} or simply an \textit{elliptic progression of primes}, if $E/\Q$ is understood from the discussion. If $p_1<p_2$ are two primes such that $\np{p_1}=\np{p_2}$, then I say that  $p_1<p_2$ are \textit{elliptic twin primes}. Thus, the assertion of \Cref{th:main} for $k=1$ is equivalent to the existence of arbitrarily long elliptic progressions of primes (for any given elliptic curve $E/\Q$). \Cref{con:k-one-case}{\bf(3)}, on the other hand, asserts that for any elliptic curves $E/\Q$, there exists infinitely many pairs of elliptic twin primes $p_1<p_2$ satisfying $\np{p_1}=\np{p_2}$. 

Numerical evidence for \Cref{con:k-one-case} is provided in \Cref{se:numerical-data}. Numerical examples show that an elliptic progression of primes need not be contained in an arithmetic progression of primes.

\Cref{con:k-one-case}{\bf(1)} seems difficult at the moment. Even the weakest assertion, namely \Cref{con:k-one-case}{\bf(3)}, seems difficult at the moment (even with GRH). These conjectures (\Cref{con:k-one-case}{\bf(1)--(3)}) point to new phenomena in the theory of prime numbers as they relate to the arithmetic of elliptic curves. On the other hand, as noted in \Cref{re:CM-k-1-case}, the $k=1$ case for a CM elliptic curve  is a version of the classical $k=2$ problem considered by \cite{erdos36b} and hence one expects that \Cref{con:k-one-case}{\bf(1)} to be true and provable unconditionally using \cite{wilson1969} and methods of \cite{erdos36b} and this paper.

Let me also point out that in \citep{silverman11} it has been conjectured that that there exists infinitely many pairs of primes $p,q$ such that $\np{p}=q$ and $\np{q}=p$ (I am grateful to Joseph Silverman for providing this reference). This conjecture is, of course, quite distinct from \Cref{con:k-one-case}.

I thank M.~Ram Murty for his comments which have improved the readability of this paper and for pointing out that in \citep{murty84a} and \citep{murty84b}   estimation of higher moments (used here and calculated by the sieving method of \citep{granville07}) is carried using a different method. 

It is a pleasure to thank the referee for suggesting improvements, simplifications, corrections and also for remarks on \Cref{con:k-one-case}{\bf(1)} which led me provide additional numerical data given in \Cref{se:numerical-data}.

\section{Basic Estimates}
Throughout this paper $\vare$ will be  a positive number which will be sufficiently small.  Symbols $c_1,c_2,\ldots$ will be positive constants as will $c_1(\vare),c_2(\vare),\ldots$ which will depend  on $\vare$. In either of the cases the  numerical values of these constants will be immaterial and I caution the reader that in some situations the same constant may be denoted by different symbols. The letter $x\geq 1$ will denote a real number and symbols $x_0,x_0(\vare),x_1,x_1(\vare),\ldots$ will also be real numbers whose values will not be important to us except for the fact that such numbers will exist (in the contexts where they appear) and again such symbols may occur in multiple contexts (and may differ from the ones which appearing in other contexts). Hopefully there will be no confusion caused by this notational conflation which I will indulge in through out this paper.

Let $E$ be an elliptic curve over $\Q$. Let $E[m]$ be the $m$-torsion of the elliptic curve $E$ for an integer $m\geq 1$. If $E/\Q$ has complex multiplication by $K$, then assume that $E$ has CM by the ring of integers $\O_K$ of the imaginary quadratic field $K$. The extension $\Q(E[m])/\Q$ ($m\geq 1$) will be called the \textit{$m$-division field of $E$}. 

By \cite[Appendix, Theorem 1]{cojocaru-kani2005} (non-CM case) and \cite{serre1972} (if $E$ has CM by the ring of integers $\O_K$ of $K$), there exists an integer $A_E\geq 1$ depending on $E$, such that for all integers $m\geq 1$ coprime to $A_E$, one has an isomorphism of finite groups
\be\label{eq:gal-grp}
 \gal(\Q(E[m])/\Q)\isom 
 	GL_2(\Z/m),
 \ee
and if $E/\Q$ has complex multiplication by ring of integers $\O_K$ of $K$, then
\be\label{eq:gal-grp2}
\gal(\Q(E[m])/K)\isom 
	(\O_K/(m))^*.
\ee

\newcommand{\sP}{\mathscr{P}} 

Let $\pi(x)$ be the number of primes $\leq x$. A prime $p$ will be called a \textit{good prime for} $E$, if $p\nmid A_E$ and if $E$ has good reduction at $p$. Let $\sP=\sP_E$ be the set of all the good primes of $E$.  Thus for any integer $m\geq 1$ divisible only by primes $p\in\sP$, one has the isomorphism \eqref{eq:gal-grp}.

For an integer $0<n$, let
\be\label{def:omegap} \omega_\sP(n)=\sum_{\sP\ni q|n}1\ee
be the number of distinct prime divisors $q\in \sP$ of $n$.  For a real number $z>0$, and an integer $n>0$, define 
\newcommand{\omegap}{\omega_{\sP}}
\newcommand{\omegapz}{\omega_{\sP,z}}
\be\label{def:omegapz}\omega_{\sP,z}(n) = \sum_{\sP\ni q|n, q<z}1,\ee so that $\omega_{\sP,z}(n)$  counts the number of distinct prime divisors  $q|n$ such that $\sP\ni q<z$.

Note that if $\sP$ is the set of all prime numbers, then 
\be\label{eq:usual-omega} 
\omegap(n)=\omega(n)
\ee 
is the usual number of all prime divisors of $n$ and 
\be\label{eq:usual-omega2} 
\omegapz(n)=\omega_z(n)
\ee
 is the usual number of all prime divisors of $n$ which are less than $z$.

\begin{defn}\label{def:grh}
	For any integer $n\geq 2$, let $E[n]$ be the $n$-torsion subgroup of $E$ and $K_n=\Q(E[n])$ be the field containing all the $n$-torsion. Then $K_n/\Q$ is a finite Galois extension with $\gal(K_n/\Q)\subseteq GL_2(\Z/n)$ and one says that \emph{Generalized Riemann Hypothesis (GRH) holds for $E/\Q$} if for all $n\geq 2$, the Dedekind zeta function $\zeta_{K_n}(s)$ has all its non-trivial zeros on the critical line $\Re(s)=\frac{1}{2}$. 
\end{defn}

Let $0<d$ be a positive, square-free integer. Let $\pi_E(x,d)$ be the number of good primes $p\leq x$ such that $d|\npp$. Write
\be 
\pi_E(x,d)=\delta(d)\lix+r_d,
\ee
where $\lix=\int_2^x\frac{dt}{\log(t)}$, where $\delta(d)$ is a positive fraction  (whose explicit form is not needed for the moment) and $r_d$ is the ``error term'' which is simply defined as $r_d=\pi_E(x,d)-\delta(d)\lix$.  Note that strictly speaking $r_d$ depends on $x$, but my notation for $r_d$ suppresses this dependence and follows standard conventions--for instance \cite{cojocaru05}. One is interested in this expression for $d$ squarefree and divisible only by the good primes $\sP$ of $E$.

First assume $E$ does not have $CM$. Then function $\delta(d)$ is multiplicative but not completely multiplicative and has been explicitly described, using \eqref{eq:gal-grp}, by many authors \cite{steuding05a,steuding05b}, \cite{liu05}, \cite{liu06} and \cite[Proposition 10]{cojocaru05}. From \cite[Page 345]{steuding05a} (note that my $\delta(d)$ is $\frac{1}{\delta(d)}$ in \cite{steuding05a}) or \cite[Proposition 10]{cojocaru05} that for any good prime $\ell$, one has
\be\label{eq:deltal1} 
\delta(\ell)=\frac{\ell^2-2\ell}{(\ell^2-\ell)(\ell^2-1)}=\frac{1}{\ell}+O(\frac{1}{\ell^2})
\ee
and more generally for any integer $\geq 1$ one has
\be\label{eq:deltal2} 
\delta(\ell^m)=\frac{1}{\ell^m}+O(\frac{1}{\ell^{m+1}}).
\ee
This is adequate for applying the sieving argument of \cite{granville07} in the non CM case. 

Now suppose $E/\Q$ has complex multiplication. Then by \cite[Page 346]{steuding05a}, \cite[Propostion 10(ii)]{cojocaru05} one has 
\be\label{eq:deltal-cm} 
\delta(\ell^m)=\begin{cases}
	\frac{1}{2\ell^{m+1}}+O(\frac{1}{\ell^{m+2}}) & \text{ if } \ell\text{ is inert in } K,\text{ and }m\text{ is odd},\\
	\frac{1}{2\ell^{m}}+O(\frac{1}{\ell^{m+1}}) & \text{ if } \ell\text{ is inert in } K,\text{ and }m\text{ is even},\\
	\frac{m+1}{2\ell^m}+O(\frac{1}{\ell^{m+1}}) & \text{ if } \ell\text{ splits in } K,\\
	\frac{1}{2\ell^m}+O(\frac{1}{\ell^{m+1}}) & \text{ if } \ell\text{ ramifies in } K.
\end{cases}
\ee

If we assume that GRH holds for $E$ then the precise form of the error term $r_d$ can be made explicit (this is the only point where one needs to assume GRH). This is done by means of the explicit Chebotarev density theorem of \citep{murty88,serre82} and hence one has the following estimate for $r_d$ due to \citep{steuding05a,steuding05b} and \cite[Proposition 11]{cojocaru05}:
\be\label{eq:error-term}
|r_d|=\begin{cases}
	O\left(d^{3/2}x^{1/2}\log(dx) \right) & \text{if } E \text{ does not have CM and},\\
	O\left(d^{1/2}x^{1/2}\log(dx) \right) & \text{if } E \text{ has CM}.
\end{cases}
\ee
and the implied constant is dependent on $E$ (specifically the conductor of $E$). 

In the CM case, I will not assume GRH and hence one cannot use the explicit error term provided by \eqref{eq:error-term}. One important point to note in both the cases, one does not need the explicit error term given by \eqref{eq:error-term}, but one needs bounds on sums of the form $\sum \abs{r_d}$ for $d$ in a suitable range. 

The following consequence of  \citep[Theorem 1]{liu06} (which treats both the CM and the non-CM cases respectively) will be crucial in the proof of \Cref{th:main}. It is important to note that \cite[Theorem 1]{liu05} (and \cite{liu06}) shows that the theorem given below is unconditional (i.e. without assuming GRH) if $E/\Q$ has complex multiplication. As is shown in \cite{liu05}, in the CM case, one can work with the Bombieri-Vingradov Theorem \cite[Theorem 3]{wilson1969} to arrive at reasonably good bounds on such sums (in the CM case). This is how one overcomes the absence of GRH in the CM case.

\bthm\label{th:liu}
\grh\ Then there exists a positive constant $c_i$ such that for all $x$ sufficiently large one has
$$\sump x (\omega(\npp)-\llx)^2=\cnst x\llx(1+o(1)).$$ 
\ethm

\Cref{th:liu} says that the numbers $\omega(\npp)$ follow a normal distribution with mean $\llx$ and variance $(\llx)^{1/2}$ (a more precise version of this is given by \cite[Theorem 1]{liu06}). The next lemma computes the fourth moment of this distribution. 
\bpro\label{pro:fourth-moment}\grh\ Then there is a constant $c_i$ such that 
for all $x$ sufficiently large, one has
$$
\sump x \left(\omegap(\npp)-\llx\right)^4=\cnst \pi(x)(\llx)^2+o(\pi(x)(\llx)^2).
$$
\epro
\bp 
This is proved by the sieving argument of \citep{granville07} using methods of \cite{liu06} and especially using \cite{liu05} to achieve unconditional (i.e. without GRH) result in the CM case. In \cite{liu05}, the use of GRH is mitigated by the use number field version of Bombieri-Vinogradov Theorem established in \cite{wilson1969}. An alternative method to this using \citep{murty84a} and \citep{murty84b} was suggested to me by M.~Ram Murty.  

First assume $E/\Q$ does not have complex multiplication. Let $\sP$ be the set of primes as defined earlier. 
Let me write $z=x^\beta$ where $\beta<1$ is a constant. At the end of this proof one is led to the choice  $\beta=\frac{1}{29}$. Let $\omega_z(n)$ be as defined in \eqref{def:omegapz}. Note that in the non-CM case, $\sP$ a omits a finite number of primes from the set of all prime numbers and hence  the difference 
$$\omega_z(n)-\omegapz(n)=\omega_z(n)-\sum_{\sP\ni q<z}1 =O(1)$$
is absolutely bounded independently of $n$ and hence makes no difference if we consider
$\omega_z(n)$ or $\omegapz=\sum_{\sP\ni q<z}1$. The latter type of sum is considered in \cite[Proposition 3]{granville07} instead of the former and to be absolutely pedantic, this is also the case in \cite{steuding05a}, \cite{cojocaru05}, \cite{liu06} but this point is not always emphasized in their notation and statements.

Further note that by Mertens' Theorem \cite[Mertens Theorem]{hardy-book}, and \eqref{eq:deltal1} one has 
\be \sum_{\sP\ni \ell<x}\delta(\ell)=\sum_{\sP\ni \ell<x}\frac{1}{\ell}=\llx +o(1).\ee

As $z=x^\beta$, any $n\leq x$ has at most a bounded number of prime divisors $q\geq z$. So one has
\be\label{eq:omegaz-trans}
\begin{aligned} 
\omegap(\npp)-\llx & =& (\omegapz(\npp)-\llx)
					  + (\omegap(\npp)-\omega_z(\npp)) \\
				  & = &  (\omegapz(\npp)-\llx) + O(1).
\end{aligned}
\ee
Taking fourth powers of \eqref{eq:omegaz-trans} one has
\beas 
(\omegap(\npp)-\llx)^4 & = & ((\omegapz(\npp)-\llx) + O(1))^4.
\eeas
So one has 
\newcommand{\oo}{\operatorname{O}\nolimits}
\be
\begin{aligned}
\sump x (\omegap(\npp)-\llx)^4 & =  \sump x (\omegapz(\npp)-\llx)^4+\\
			&  \qquad+ \oo\left( \sump x (\omegapz(\npp)-\llx)^3 \right)+ \\ 
			&   \qquad+ \oo\left( \sump x (\omegapz(\npp)-\llx)^2 \right)+ \\
			&   \qquad+ \oo\left( \sump x (\omegapz(\npp)-\llx)^1 \right)+ \\
			&   \qquad+ \oo\left(\pi(x)\right). 					
\end{aligned}
\ee
Each of the terms on the right will be estimated using the sieving method of \cite{granville07}. For $i=1,2,3,4$, let
\be \label{def:d4} 
D_i= \left\{d=q_1\cdot q_2\cdots q_i: \text{ with } \sP\ni q_j<z \text{ for } j=1,\cdots, i \right\}
\ee be the set of square numbers which are products of at most four distinct primes (from $\sP$) $q_1,\ldots, q_4<z$.
By \citep[Proposition 3]{granville07} one has
\be
\begin{aligned}
\sump x (\omegapz(\npp)-\llx)^4 &= \cnst \pi(x)(\llx)^2\left(1+\oo\left(\frac{1}{\llx}\right)\right)\\
                               &\qquad +\oo\left((\llx)^4\sum_{d\in D_4} |r_d|\right),
\end{aligned}
\ee
where the set $D_4$ is defined by \eqref{def:d4} for $i=4$. Similarly, by \cite[Proposition 3]{granville07} one has
\be
\begin{aligned} 
\sump x (\omegapz(\npp)-\llx)^3 &\leq \cnst \pi(x)\llx\left(1+\oo\left(\frac{1}{\llx}\right)\right)\\
&\qquad +\oo\left((\llx)^3\sum_{d\in D_3} |r_d|\right),
\end{aligned}
\ee
where  $D_3$ is the set defined by \eqref{def:d4} for $i=3$. 

The sum \be \sump x (\omegapz(\npp)-\llx)^2\ee is estimated by \Cref{th:liu} and is $\oo(\pi(x)\llx)$ and so by Cauchy-Schwarz inequality one has 
\be\bigg\vert\sump x (\omegapz(\npp)-\llx)\bigg\vert\ll \pi(x)\llx.\ee
Let me indicate how to estimate the term (the estimation of other terms is similar to this one)
\be\label{eq:error-erdos-kac} 
\oo\left((\llx)^4\sum_{d\in D_4}|r_d|\right)
\ee  
where the sum over $d$ runs over all square-free integers which are products of at most four primes $\sP\ni q<z$.  Each $d\in D_4$ which contributes to the sum is at most a product of four primes $\sP\ni q<z$ so
\be 
d=q_1\cdots q_4 <z^4,
\ee
and hence $d^{3/2}\log(d)\leq d^{5/2}<(z^4)^{(5/2)}=z^{10}$. Further since there are $\pi(z)$ primes less that $z$, so the number of such $d$ is $O(\pi(z)^4)$.  Thus, one has
\be 
\sum_{d\in D_4}|r_d| = O\left(x^{1/2}\lx z^{10}\pi(z)^4 \right)=O\left(x^{1/2}\lx \frac{z^{14}}{(\log(z))^4}\right)
\ee 
which is
\be 
O\left(\frac{x^{1/2} x^{14\beta}}{(\lx)^3}\right)
\ee
Now choose $\beta$ so that $14\beta+\frac{1}{2}<1$ which gives $\beta<\frac{1}{28}$. So if we choose $\beta=\frac{1}{29}$ one  gets
\be 
\sum_{d\in D_4}|r_d|=O\left(\frac{x^{57/58}}{(\lx)^3}\right),
\ee
and finally with this estimate the last term in \eqref{eq:error-erdos-kac} is
\be 
\leq \cnst \frac{x^{57/58}(\llx)^4}{(\lx)^3} = o(\pi(x)(\llx)^2).
\ee
The sum
\be\label{eq:error-erdos-kac2} 
\oo\left((\llx)^3\sum_{d\in D_3}|r_d|\right)
\ee
can also be estimated similarly and is $o(\pi(x)(\llx)^2)$. The sum \be\oo\left(\llx\sum_{d\in D_1}|r_d|\right)\ee can also be estimated similarly and is certainly $o(\pi(x)(\llx)^2)$. This completes the proof of \Cref{pro:fourth-moment}.

Now let me sketch a proof in the case $E/\Q$ has CM. As was noted earlier, the set $\sP$ in the CM case consists of primes $p$ such that $E$ has good reduction at $p$, and $p$ splits in the CM field $K$ of $E$. This means that $p$ is a prime of good ordinary reduction for $E$. Thus in the CM case, one is restricting oneself to the case of primes of good ordinary reduction. As is well known, the set of such rational primes  have density one half in all the rational primes. The choice of this set $\sP$ is permitted in the framework of \cite[Proposition 3]{granville07} and its proof. By \eqref{eq:deltal-cm} one knows that 
$$\sum_{\sP\ni \ell\leq x}\frac{1}{\ell}=\frac{1}{2} \llx+o(1).$$
As is remarked in \cite{granville07}, and shown in \cite{liu05}, the absence of GRH in the CM case is circumvented by means of the Bombieri-Vinogradov Theorem for number fields  of \cite{wilson1969} (applied here for the CM field $K$). Hence one can apply \cite[Proposition 3]{granville07}, and estimate the error term sums using  \cite{liu05} and proceed exactly as above to arrive at the claimed assertion in the CM case. 
\ep
 
The following theorem was proved in \citep{erdos36} for $k=3$. 
\bthm\label{th:erdos-product}
Let $k\geq 3$ be an integer. Let $0<\vare < \frac{1}{2}-\frac{1}{k\log 2}$ be a  real number. Let
\be\label{eq:k-fact} 
\dx=\left\{n=n_1\cdots n_k\leq x: \omega(n_i)>\low \ {\rm for\ } 1\leq i\leq k \right\}.
\ee
Then 
\be 
\cdx \leq \cnste \frac{x}{(\lx)^{1+\delta}},
\ee
with
$\delta(k)=(\log 2)k(1-2\vare)-2>0$ for all $k\geq 3$.
\ethm
\bp 
Let $n\in \dx$. Let us write 
\be\label{eq:m1m2-fact} n=m_1\cdot m_2\ee
where $m_1$ is square-free and for every prime $q|m_2$, one has $q^2|m_2$. Let 
\be\label{eq:dxo-def}\dxo=\left\{n\in \dx: \omega(m_2)>\vare\llx \right\}\ee
and let 
\be \label{eq:dxt-def}\dxt=\left\{n\in \dx: \omega(m_2)\leq\vare\llx \right\}.\ee
Then clearly $\dx=\dxo\cup\dxt$ and $\dxo\cap\dxt=\emptyset$. The cardinalities of each of $\dxo$ and $\dxt$ will be estimated separately.

If $n\in\dxo$ then there is a divisor $d|n$ such that $\omega(d)=j$ where $j=[\vare\llx]$ and for each prime $q|d$ one has $q^2|n$. Hence,  $\cdxo$ is less than the number of such numbers less than $x$ i.e., 
\be\label{eq:K0} 
\cdxo <\cnst x\left(\sum_{q,\alpha\geq 2}\frac{1}{q^\alpha}\right)^j\frac{1}{j!},
\ee 
where the sum is over all primes $q$ and all integers $\alpha\geq 2$ and clearly the sum converges so let 
$K=\sum_{q,\alpha\geq 2}\frac{1}{q^\alpha}$. Then \eqref{eq:K0} is
\be\label{eq:j-K-eq}  
\cdxo <\cnst \frac{x K^j}{j!}.\ee

For any integer $j\geq 1$, define $a_j=K^j/j!$. Then $$\frac{a_{j+1}}{a_j}=\frac{K^{j+1}}{(j+1)!}/\frac{K^j}{j!}=K/(j+1),$$
and 
$$\frac{a_{j+2}}{a_{j+1}} = \frac{K}{j+2} <\frac{K}{j+1}=\frac{a_{j+1}}{a_j}.$$
Thus, $a_{j+1}/a_j$ is a positive, monotone decreasing sequence. Let $\delta_1>1$, then for all
sufficiently large positive $j$ one has
$\frac{a_{j+1}}{a_j}\leq \frac{1}{e^{\delta_1}}$. Thus $a_{j+1}=e^{-\delta_1}\cdot a_j$ and thus by induction on $j$ one see that $a_j\leq e^{\delta_1\cdot j}$. As $j=\vare \llx$ one obtains that 
\be\label{eq:j-K-eq2} a_j=O(e^{-\delta_1\cdot\vare \cdot \llx})=O((\lx)^{-\delta_1\vare}).\ee
Hence, combining \eqref{eq:j-K-eq2} with \eqref{eq:j-K-eq} one obtains that
\be  
\cdxo <\cnst \frac{x K^j}{j!}<\cnst \frac{x}{(\lx)^{\delta_1\vare}}
\ee
and where for any fixed $\varepsilon$, one can choose $\delta_1>1$ so  that the equation $\delta_1\varepsilon=1+\delta>1$ has a solution for any arbitrary positive $\delta$.
Hence one deduces that for all $x$ sufficiently large, one has
\be\label{eq:cdxo-estmiate} 
\cdxo \leq \cnst x \frac{K^j}{j^j} \leq \cnst \frac{x}{(\lx)^{1+\delta}}
\ee
for any positive $\delta$.

Now let me estimate $\cdxt$. By definition \eqref{eq:k-fact}, each $n\in \dxt$ is of the form  $n=n_1\cdots n_k$ with $\omega(n_i)\geq \low$ for all $1\leq i\leq k$. On the other hand from \eqref{eq:m1m2-fact} one knows that
$n=m_1m_2$ where $m_1$ is square-free and for every $q|m_2$, $q^2|m_2$. Hence, I claim that
\be\label{eq:omega-bnd}  
\omega(n)\geq \omega(n_1)+\omega(n_2)+\cdots+\omega(n_k)-(k-1)\omega(m_2).
\ee
This is clear from the two ways of writing $n$ given by \eqref{eq:k-fact} and \eqref{eq:m1m2-fact} and the fact that if a prime divides two or more of the $n_1,\ldots,n_k$ then it counts at most $k$ times in the sum  $\omega(n_1)+\omega(n_2)+\cdots+\omega(n_k)$ and it is also counted in $\omega(m_2)$ as it divides $m_2$.

Thus, from \eqref{eq:dxt-def} and \eqref{eq:omega-bnd} one obtains
\be
\begin{aligned}
\omega(n) & \geq k(1-\varepsilon)\llx-(k-1) \varepsilon\llx.\\
			& \geq k(1-2\varepsilon)\llx+\varepsilon\llx\\
			& \geq k(1-2\varepsilon)\llx.
\end{aligned}
\ee
The number on the right is positive if $\varepsilon<1/2$.
Therefore, it follows that if $\varepsilon<1/2$, then each $n\in\dxt$ has at least $k(1-2\vare)\llx>0$ distinct prime factors. Thus, each such $n$ has at least $2^{k(1-2\vare)\llx}$  distinct divisors i.e.
\be  
d(n)\geq 2^{k(1-2\vare)\llx}.
\ee

On the other hand, from \cite[Chap. XVIII, Theorem 320]{hardy-book},  one has 
\be 
\sum_{n\leq x}d(n) \leq\cnst x\lx.
\ee
Hence, one has
\be 
2^{k (1-2\vare)\llx}\sum_{n\in\dxt} 1 \leq \sum_{n\leq x}d(n) \leq\cnst x\lx
\ee
which, using $2^\alpha=e^{\alpha \log 2}$, gives
\be
\begin{aligned}
\cdxt & \leq  \cnst \frac{x\lx}{2^{k(1-2\vare )\llx}} \\
      & \leq  \cnst \frac{x\lx}{e^{(\log 2)k(1-2\vare)\llx}} \\
      & \leq  \cnst \frac{x\lx}{(\lx)^{(\log 2)k(1-2\vare)}} \\
      & \leq  \cnst \frac{x}{(\lx)^{(\log 2)k(1-2\vare)-1}}.
\end{aligned}
\ee
This will be $\ll \frac{x}{(\lx)^{1+\delta}}$ with $\delta>0$ if 
\be 
(\log 2)k(1-2\vare)-1 > 1.
\ee
This equation is  satisfied provided one has
\be 
0<\vare < \frac{1}{2}-\frac{1}{k\log 2}.
\ee
If $k\geq3$, then this evidently holds true. But if $k=2$, then the right-hand side is negative as $\frac{1}{2\log2}=0.7213>1/2$ so there is no positive $\varepsilon$ which works for $k=2$. This is the only place where $k\geq 3$ is used in this proof. This is why the assumption $k\geq 3$ is necessary in \Cref{th:erdos-product}. 

Hence, in this case one has 
\be\label{eq:cdxt-estimate} 
\cdxt\leq\cnste \frac{x}{(\lx)^{1+\delta(k)}}
\ee 
where
\be\label{def:deltak} \delta(k)=(\log 2)k(1-2\vare)-2>0 \quad  \text{ for all }k\geq 3.
\ee  

Now one has $\cdx=\cdxo+\cdxt$. Putting the estimates \eqref{eq:cdxo-estmiate} and \eqref{eq:cdxt-estimate} together one obtains that for each $k\geq 3$, there is a $\delta=\delta(k)>0$ given by \eqref{def:deltak} such that for all sufficiently large $x$ one has
\be 
\cdx \leq \cnste\frac{x}{(\lx)^{1+\delta}}.
\ee
This proves the theorem.
\ep

\brem 
In \cite{erdos36b}, Erd{\H o}s uses a different method to deal with $k=2$ in his case.
\erem

\section{Proof of \Cref{th:main}}
\blem\label{le:basic-count1}
\grh\ Let $\vare$ be a fixed, sufficiently small positive real number. 
Let $\ax$ denote the number of primes $\sP\ni p\leq x$ such that 
\be 
\omegap(\npp)<\low \ee \
Then there is a positive constant $c_i(\vare)$ such that for all $x\geq x_0$ one has 
$$ 
\cax\leq\cnste \frac{\pi(x)}{(\llx)^2}.
$$
\elem
\bp 
It is sufficient to prove that the number of $p\leq x$ such that 
\be \omegap(\npp)<\low \ee is 
$\ll\frac{\pi(x)}{(\llx)^2}$ as the number of $p\leq x$ counted in the assertion of this Lemma is certainly less than the number of primes $p\leq x$   with this  property.
For $p\leq x$  satisfying the above property one has $$\omegap(\npp)-\llx < -\vare\llx$$ i.e. 
$$-(\omegap(\npp)-\llx)>\vare\llx>0$$ and hence $$(\omegap(\npp)-\llx)^4 > \vare(\llx)^4.$$ Hence
\be  
\sum_{p\in \ax}\vare^4(\llx)^4<\sump x(\omegap(\npp)-\llx)^4 \leq \cnst\pi(x)(\llx)^2,
\ee 
by \Cref{pro:fourth-moment}. Thus the assertion follows.
\ep

\blem\label{le:convergence}
\grh\ The series $\sum_{p\in\ax}\frac{1}{p}$ converges.
\elem
\bp 
This follows from the previous lemma by Abel summation formula \cite{apostol-analytic-book} and the tautological bound ${\mathcal A}(t)\leq \pi(t)$ for all $t\geq 1$ as the integral $$\int_5^x\frac{dt}{t\log(t)(\log\log(t))^2}$$ converges as $x\to\infty$.
\ep

\begin{defn}\label{def:admissible}
Let $E/\Q$ be an elliptic curve over $\Q$. Let $\varepsilon>0$ be a small positive real number.
I say that a prime $p\in\sP$ is an \textit{admissible prime for $E$} or simply an \textit{admissible prime}  if $\omegap(\np{p})\geq\low$, otherwise I say that \textit{$p\in\sP$ is an inadmissible prime}. 
\end{defn}

Note that this definition of admissible primes depends on $E$ and $\varepsilon$.

\blem\label{le:basic-count2}
\grh\ Let $\vare$ be a fixed, sufficiently small positive real number. 
Let $\cx$ denote the number of admissible primes $p\leq x$ i.e. the set of primes $\sP\ni p\leq x$ such that 
\be 
 \omegap(\npp) \geq \low.
\ee \
Then there is a positive constant $c_i(\vare)$ such that for all $x\geq x_0$ one has 
$$ 
\ccx\geq \pi(x)-\cnste\frac{\pi(x)}{(\llx)^2}>\cnst(\vare)\pi(x).
$$
\elem
\bp 
By \Cref{le:basic-count1}, the number of primes $p\leq x$ such that $ \omegap(\npp) \geq \low$ does not hold is $o(\pi(x))$ hence the assertion follows.
\ep

\blem\label{le:series-lemma}
\grh\ Suppose $\vare$ is a sufficiently small positive real number. Let $k\geq 1$ be an integer and let $b>a\cdot e$ with $0<a<b<1$  be  fixed positive real numbers (hence $\log(b)-\log(a)>1$). Let
$$ 
\sum\nolimits_1= \sum_{x^a\leq p_1,\ldots,p_{k} <x^b} \frac{1}{p_1\cdots p_{k}},
$$
where the sum is over all the pairwise distinct primes $p_1,\ldots, p_k$ such that $x^a<p_i<x^b$ and  $\omegap(\np{p_i})>\low$ for all $1\leq i\leq k$. Then there is a positive constant (depending on $a,b$) such that
$$\sum\nolimits_1 = \cnst(a,b,k)-\epsilon(x)>1,$$
where the right-hand side is greater than one for all sufficiently large $x$ and $\epsilon(x)\to0$ as $x\to\infty$. 
\elem
\bp Let me first remark that the proof shows that one may take $a,b$ satisfying
$0<b<\frac{1}{k-1}$ and $a<\frac{1}{e(k-1)}$ then one has $0<a<b<1$ and $\log(b/a)>1$. For $k=3$, in \citep{erdos36} $b=\frac{1}{4}<\frac{1}{2}$ and $a=\frac{1}{8}<\frac{1}{2e}$ so that $\log(b/a)=\log(4)>1$. Consider (the sum over all primes)
$$\sum_{x^a\leq p\leq x^b}\frac{1}{p}=\log\log(x^b)-\log\log(x^a)+o(1)=\log(b)-\log(a)+o(1)$$ 
by \citep[Mertens Theorem]{hardy-book} and in particular, by the choice of $a,b$, this sum is greater than a positive constant for all $x$ sufficiently large. 

Note that in case $E$ has CM, the set of primes $\sP$ has density one half and so the above equation needs to be modified with a factor of $1/2$ appearing on the right. This does not affect the rest of the proof in anyway--notably the assertion of \Cref{le:series-lemma} is true regardless of the two distinct choices for $\sP$ in the non-CM and CM cases.

Now write this sum as
\be  
\sum_{x^a\leq p\leq x^b}\frac{1}{p}=\sumtwo\frac{1}{p}+\sum\nolimits_3\frac{1}{p}
\ee
where 
\be 
\sum\nolimits_2=\sum\nolimits_2 \frac{1}{p}=\sum_{x^a\leq p\leq x^b,\ \omegap(\npp)> \low}\frac{1}{p}
\ee
is a sum over admissible primes \eqref{def:admissible}, i.e. over primes $p$ for which $\omegap(\npp)>\low$, and
\be 
\sum\nolimits_3=\sum\nolimits_3\frac{1}{p}=\sum_{x^a\leq p\leq x^b,\ \omegap(\npp)\leq \low}\frac{1}{p}
\ee
is a sum over inadmissible primes. By the convergence of the sum over inadmissible primes given by \Cref{le:convergence} one sees that if we write
\be 
\sum\nolimits_3=\sum_{x^a\leq p\leq x^b,\omegap(\npp)\leq \low}\frac{1}{p}=\epsilon(x)
\ee
then $\epsilon(x)=o(1)$ as $x\to\infty$ and in particular for all sufficiently large $x$ one has $|\epsilon(x)|<\frac{1}{100}$. The choice of $\frac{1}{100}$ is arbitrary, the point being that $\epsilon(x)$ can be made smaller than any given positive real number by choosing $x$ sufficiently large. So one gets that the sum of the reciprocals of admissible primes (on the left)
\be 
	\sumtwo=\sum_{x^a\leq p\leq x^b}\frac{1}{p}-\epsilon(x)=\log(b/a)-\epsilon(x)
\ee
and hence for sufficiently large $x$ the right hand side is greater than one as $\epsilon(x)$ can be made arbitrarily small by choosing $x$ sufficiently large. Thus the result is true for $k=1$.

Now before I give the general case let me illustrate the method of proof with the case $k=2$ as this also will allow me to set up notational conventions needed to do the general case. \textit{The sums over primes which occur below are sums over admissible primes.}
Since the result is true for $k=1$, one gets by squaring 
\be 
\left(\sumtwo\right)^2=
(\log(b/a)-\epsilon(x))^2
\ee
on the other hand by squaring the sum over admissible primes on the left one gets a sum over admissible primes
\be 
\left(\sumtwo\right)^2=\sum \frac{1}{p_1p_2}+2\sum \frac{1}{p_1^2}=(\log(b/a))^2-2\log(b/a)\epsilon(x)+\epsilon^2(x).
\ee
Now the sum $\sum_p\frac{1}{p^2}$, taken over all primes (and not just admissible ones) is convergent and hence the sum $\sum_{x^a<p_1<x^b}\frac{1}{p_1^2}$ can be made arbitrarily small by choosing $x$ sufficiently large, in other words the sum is $\epsilon(x)$ for some function $\epsilon(x)$.  Since $\log(b/a)$ is fixed one may also write middle term in the above expression as $\epsilon(x)$ for yet another $\epsilon(x)$. Since there are finitely many such $\epsilon(x)$ it is possible to arrange $x$ so large that all of them be less than a given positive real number, similarly $\epsilon^2(x)$ is also $\epsilon(x)$. \textit{Note that I am systematically conflating all the $\epsilon(x)$--a diligent reader may easily workout the full notational genuflections required to separate them.}
Thus one gets from all this that
\be 
\sum_{x^a<p_1,p_2<x^b}\frac{1}{p_1p_2}=\cnst-\epsilon(x)>1
\ee
for all sufficiently large $x$.

Now let me return to the general case which is proved by induction on $k$. Suppose the result is true for such a sum in the statement of the lemma and all natural numbers $\leq k$. Then I show that the result is also true for $k+1$. Since the result is true for $k$ one has
\be 
\sum\nolimits_1= \sum_{x^a\leq p_1,\ldots,p_{k} <x^b} \frac{1}{p_1\cdots p_{k}}=\cnst(a,b,k)-\epsilon(x),
\ee
Then consider the product
\be \label{eq:product-sum}
\left(\sum_{x^a\leq p_1,\ldots,p_{k} <x^b} \frac{1}{p_1\cdots p_{k}}\right)
\left(\sum_{x^a\leq p_1,\ldots,p_{k} <x^b} \frac{1}{p}\right)
\ee
where both the sums are over all (pairwise) distinct admissible primes in the asserted range. Then by induction hypothesis (applied to both the sums) one has
\be\label{eq:ind-step} 
\left(\sum_{x^a\leq p_1,\ldots,p_{k} <x^b} \frac{1}{p_1\cdots p_{k}}\right)
\left(\sum_{x^a\leq p_1,\ldots,p_{k} <x^b} \frac{1}{p}\right) = (\cnst-\epsilon(x))(\cnst-\epsilon(x))=\cnst-\epsilon(x),
\ee
for some $\epsilon(x)$  and some constant depending on $a,b$ and in particular the right hand side is greater than one.

Now multiplying out the product one gets one term of the form
\be 
\sum\nolimits_1=\left(\sum_{x^a\leq p_1,\ldots,p_{k+1} <x^b} \frac{1}{p_1\cdots p_{k+1}}\right)
\ee
a sum over pairwise distinct admissible primes in the asserted range
and the remaining terms are of the form
\be\label{eq:coincident-sums} 
\sum\nolimits_1=\left(\sum_{x^a\leq p_1,\ldots,p_{k+1} <x^b} \frac{1}{p_{i_1}^{\ell_1}\cdots p_{i_j}^{\ell_j}}\right),
\ee
where  the primes are admissible and $\ell_1+\cdots +\ell_j=k+1$ and $j\leq k$ and at least one of the exponents $\ell_1,\ldots,\ell_j$ is $\geq 2$, and which arise out of coincidences among the primes when one takes the product in \eqref{eq:product-sum}. So one has to deal with these two sort of sums. 

I claim that the sums \eqref{eq:coincident-sums} where one has coincidences between some of the primes, are small i.e. I claim that for any $k$ and $\ell_1+\cdots+\ell_j=k$ and at least one $\ell_i\geq 2$,
\be\label{le:subclaim} 
\sum\nolimits_1=\left(\sum_{x^a\leq p_1,\ldots,p_{k+1} <x^b} \frac{1}{p_{i_1}^{\ell_1}\cdots p_{i_j}^{\ell_j}}\right)=\epsilon(x),
\ee
where $\epsilon(x)\to0$ as $x\to\infty$.

Suppose the assertion \eqref{le:subclaim} is true for the moment. Then one gets from \eqref{eq:ind-step} that
\be\label{eq:subsubclaim1}  
\left(\sum_{x^a\leq p_1,\ldots,p_{k} <x^b} \frac{1}{p_1\cdots p_{k}}\right)
\left(\sum_{x^a\leq p_1,\ldots,p_{k} <x^b} \frac{1}{p}\right)  =  \cnst-\epsilon(x).
\ee
The left-hand side of \eqref{eq:subsubclaim1} can be written as
\be\label{eq:subsubclaim2}
\begin{aligned}
\left(\sum_{x^a\leq p_1,\ldots,p_{k} <x^b} \frac{1}{p_1\cdots p_{k}}\right)
\left(\sum_{x^a\leq p_1,\ldots,p_{k} <x^b} \frac{1}{p}\right) & =  \\ &= \left(\sum_{x^a\leq p_1,\ldots,p_{k+1} <x^b} \frac{1}{p_1\cdots p_{k+1}}\right) + \\
&\qquad + \sum_{j,\ell_1,\ldots,\ell_j} \left(\sum_{x^a\leq p_{i_1},\ldots,p_{i_j} <x^b} \frac{1}{p_{i_1}^{\ell_1}\cdots p_{i_j}^{\ell_j}}\right)\\
&= \left(\sum_{x^a\leq p_1,\ldots,p_{k+1} <x^b} \frac{1}{p_1\cdots p_{k+1}}\right)+\epsilon(x),
\end{aligned}
\ee
and hence rearranging terms of \eqref{eq:subsubclaim1} and \eqref{eq:subsubclaim2}, one has
\be 
\left(\sum_{x^a\leq p_1,\ldots,p_{k+1} <x^b} \frac{1}{p_1\cdots p_{k+1}}\right)=\cnst -\epsilon(x)>1.
\ee
Thus, to complete the proof of \Cref{le:series-lemma}, it remains to prove  the claim \eqref{le:subclaim} regarding the sums \eqref{eq:coincident-sums}. This claim is again proved by induction on $j$. Suppose $j=1$. Then one has $\ell_1\geq 2$ as at least one of the exponents is required to be greater than one. Then the assertion follows as $\sum\frac{1}{p^2}$ converges. Hence assume the result is true when the number of primes factors is less than or equal to some $j$ and any $\ell_1,\ldots,\ell_j$ with one of them being greater that two. Then I prove the result is true for $j+1$ and any $\ell_1,\ldots,\ell_{j+1}$ with at least one of them being greater than two.  

If at least one of $\ell_i\geq 2$ for some $1\leq i\leq j$ then the assertion follows by induction as $\sum_{x^a<p\leq x^b}\frac{1}{p^{\ell_{j+1}}}$ is bounded if $\ell_{j+1}=1$ and $o(1)$ if $\ell_{j+1}\geq 2$. 

So now assume that $\ell_{i}=1$ for all $1\leq i\leq j$ and $\ell_{j+1}\geq 2$. Then the sum in consideration  occurs in the product 
\be
\left(\sum_{x^a<p_i\leq x^b}\frac{1}{p_1\cdots p_j}\right)\left(\sum_{p}\frac{1}{p^{\ell_{j+1}}}\right)
\ee
and the first is a summand of $\left(\sum_{x^a<p\leq x^b}\frac{1}{p}\right)^j$ and hence is bounded (again one uses induction hypothesis to deal with the series with coincidences as the number of factors will be $<j+1$) and the second sum is $o(1)$ as $\sum_p\frac{1}{p^{\ell_{j+1}}}$ converges as $\ell_{j+1}\geq 2$. Thus, in all cases the claim \Cref{le:subclaim} is proved and hence \Cref{le:series-lemma} is proved. 
\ep

\bp[{Proof of \Cref{th:main}}]
\grh\ Let $\vare$ be a fixed, sufficiently small positive real number. Let $k\geq 1$ be an integer.  Let 
\be 
\bx{k} =\{  p_1\cdot p_2\cdots p_k\leq x:   p_i\neq p_j\ {\rm if}\ i\neq j\  \ {\rm and}\ p_i\ {\rm admissible\ for\ all}\    1\leq i\leq k\}.
\ee
Then I claim that it is sufficient to prove that for any $k\geq 1$ one has
\be\label{eq:erdos-bound}
\cbx{k}\geq \cnste\frac{x}{\lx}.
\ee
Note that this claim does not use the assumption in \Cref{th:main} that $k\geq 3$. This assumption enters the proof only through the use of \Cref{th:erdos-product} which needs $k\geq 3$.

Suppose for the moment that \eqref{eq:erdos-bound} has been established. Then one can complete the proof of the theorem as follows. On one hand, one has $\cbx{k}\geq \cnste \frac{x}{\lx}$ and on the other hand, by \Cref{th:erdos-product}, one has $\cdx\leq \cnste\frac{x}{(\lx)^{1+\delta}}$. 

So if one considers the mapping 
\be\label{eq:mapping}
\bx{k}\to \dx\ee 
given by 
\be p_1\cdots p_k\in \bx{k} \mapsto \np{p_1}\cdots \np{p_k}.\ee 
Then by definition, the number $G_k(n)$ is precisely the cardinality of the fiber over $n$ under the mapping given by \eqref{eq:mapping}. By \eqref{eq:erdos-bound} and \Cref{th:erdos-product},  some fibers of the mapping \eqref{eq:mapping} must have cardinality $\gg (\log x)^\delta$. Otherwise, for any arbitrary real number $c>0$, the cardinality of any fiber is $\leq c\cdot(\log x)^\delta$. This together with
\Cref{th:erdos-product}, gives
\be \cbx{k}\ll \cdx\cdot c\cdot (\log x)^\delta \leq c\cdot \frac{x}{\log(x)},\ee 
which contradicts \eqref{eq:erdos-bound} because $c>0$ is arbitrary. Hence, there exist $1\leq n\leq x$ such $G_k(n)\gg (\log x)^\delta$. 

So let me now prove the bound asserted in \eqref{eq:erdos-bound}. If $k=1$ then this follows from \Cref{le:basic-count1} and \Cref{le:basic-count2}. So assume $k\geq 2$. If $m=p_1\cdots p_k\leq x$ then $p_k=\frac{m}{p_1\cdots p_{k-1}}\leq \frac{x}{p_1\cdots p_{k-1}}$ and hence $p_k\in \bxx{1}{x/(p_1\cdots p_{k-1})}$. Conversely if $p_1,\ldots,p_{k-1}$ are admissible primes and $p_1\cdots p_{k-1}\leq x$ and $p_k\leq \frac{x}{p_1\cdots p_{k-1}}$ is an admissible prime then $p_1\cdots p_k\leq x$ and is a product of admissible primes and hence is in $\bx{k}$. Moreover any of the $k!$ permutations of the factors $p_1,\ldots, p_k$ give the same $m=p_1\cdots p_k$. So one has
\be 
k! \cbx{k} \geq \sum \cbxx{1}{x/(p_1\cdots p_{k-1})},
\ee
where the sum is over all the pairwise distinct $p_1,\ldots,p_{k-1}$ and excludes the terms corresponding to the coincidences $p_k=p_i$ for any of the $1\leq i\leq k-1$. Since one wants a lower bound, it is in fact enough to replace this sum by one of the same sort except where the primes $p_i$ are in some fixed range $x^a\leq p_i <x^b$ for $0<a<b<1$, and which satisfy the requirements of \Cref{le:series-lemma}.
Thus
\be 
k! \cbx{k} \geq \sum_{x^a\leq p_1,\ldots,p_{k-1} <x^b} \cbxx{1}{x/(p_1\cdots p_{k-1})},
\ee

If one chooses $b$  so that  
\be\label{eq:b-bound} 
1-b(k-1)>0
\ee
 then $x/(p_1\cdots p_{k-1})\geq x/x^{b(k-1)}\gg 1$ 
and one can apply \Cref{le:basic-count2} to get
\be 
k! \cbx{k} \geq \sum_{x^a\leq p_1,\ldots,p_{k-1} <x^b} \cbxx{1}{x/(p_1\cdots p_{k-1})}\geq \cnste \sum_{x^a\leq p_1,\ldots,p_{k-1} <x^b} \pi(x/(p_1\cdots p_{k-1})).
\ee 
Hence
\be 
k! \cbx{k} \geq \cnste\sum_{x^a\leq p_1,\ldots,p_{k-1} <x^b} \frac{x/(p_1\cdots p_{k-1})}{\log(x/(p_1\cdots p_{k-1}))}.
\ee
Since $\frac{1}{\log(x/(p_1\cdots p_{k-1}))}\geq \frac{1}{\lx}$ this gives
\be 
k! \cbx{k} \geq \cnste\frac{x}{\lx}\left(\sum_{x^a\leq p_1,\ldots,p_{k-1} <x^b} \frac{1}{p_1\cdots p_{k-1}}\right).
\ee
Now the proposition follows from \Cref{le:series-lemma} provided that in addition to the condition \eqref{eq:b-bound} one choose $a$ so that the conditions of \Cref{le:series-lemma} are met. Thus on choosing $0<b<\frac{1}{k-1}$ and $0<a<\frac{1}{e(k-1)}$ one gets the assertion of  \Cref{th:main}.
\ep

In case $E$ has CM, the proof of \Cref{th:main} shows a little more. Let $G_k^{ord}(E,n)$ be the number of solutions  to  \eqref{eq:erdos} in primes of good ordinary reduction for $E$. Then one has proved the following:

\bcor\label{cor:cm-case} 
Let $E/\Q$ be an elliptic curve with complex multiplication. Then there exists integers $n\to\infty$ such that the number of solutions, $G_k^{ord}(E,n)$ to  \eqref{eq:erdos} in primes of good ordinary reduction also goes to infinity i.e. 
$$ 
\varlimsup_{n\to\infty} G_k^{ord}(E,n)=\infty.
$$
\ecor
\bp 
This is immediate from the proof of \Cref{th:main} because under the CM hypothesis, every prime $p\in\sP$ is a prime of good ordinary reduction for $E$.
\ep

\section{The cases $k=1,2$}\label{se:conj}
Let me dispense with $k=2$ case first. The $k=2$ case is robustly supported by numerical evidence. I expect that just as as in the classical case \cite{erdos36b}, the case $k=2$ requires a different treatment from the cases $k\geq 3$. Especially I expect that \Cref{th:main} holds for $k=2$ (pretty much under the same hypothesis).  This investigation is left to the readers using this paper and \cite{erdos36b} as a guide.

Now let us consider the case $k=1$. This is a case which does not occur in the classical context considered in \citep{erdos36}, \citep{erdos36b} and is even more interesting than $k\geq 2$. Let $G_1(E,n)$ be the number of solutions (in primes $p$) to the equation $N_p(E)=n$. For simplicity, write $G_1(n)=G_1(E,n)$. Then example for the curve $E:y^2+y=x^3-x$ one has $$\np{1009}=\np{1063}=1057.$$ So $G_1(1057)\geq 2$ and a simple calculation using \eqref{eq:hw-bound} and a computer search shows that equality holds so $G_1(1057)=2$.

I conjecture that for any elliptic curve $E/\Q$, one has the following: 
\begin{conj}\label{con:k-one-case}
Let $E/\Q$ be an elliptic curve and for any integer $n\geq1$ let $G_1(n)$ be the number of good primes $p$ such that $\np{p}=n$. Then
\benum
\item $\varlimsup_{n\to\infty}G_1(n)=\infty$.
\item $\varlimsup_{n\to\infty}G_1(n)\geq 2$
\item There exists infinitely many $n$ such that $G_1(n)=2$.
\eenum
\end{conj}
Clearly \Cref{con:k-one-case}(1)$\implies$\Cref{con:k-one-case}(2) and \Cref{con:k-one-case}(3)$\implies$\Cref{con:k-one-case}(2). Note that \Cref{con:k-one-case}(3) says that for any elliptic curve $E/\Q$, there are infinitely many distinct good primes $p,q$ such that $\np{p}=\np{q}$. 

In \Cref{se:numerical-data}, I provide the numerical data I have computed supporting these conjectures.

\brem 
\Cref{con:k-one-case}(1) means that for any elliptic curve $E/\Q$, there exists integers $n$ (going to infinity) and arbitrarily large strings of primes $p_1<p_2<\cdots <p_m$ such that 
\be\label{eq:prog-seq} \np{p_1}=\np{p_2}=\cdots=\np{p_m}=n.\ee
This seems substantially more difficult at this point (even on GRH), and of course it would even more astounding if this turns out to be false. It seems natural to call a sequence of primes occurring in \eqref{eq:prog-seq} \textit{an elliptic progression  of primes given by $E$} or simply \textit{primes in an elliptic progression} if the curve is unambiguously defined. If $m=2$ (resp. $m=3$), one calls such primes \textit{elliptic twin primes} (resp. \textit{elliptic prime triplets}) and so on.  At the very least, I expect that elliptic twin primes are infinite i.e. $G_1(E,n)\geq 2$ infinitely often \textit{i.e.,} there are infinitely many pairs of primes in an elliptic progression. 
\erem

The following lemma illustrates the role which primes of ordinary reduction for $E$ play in $G_1(n)$ and hence in \Cref{con:k-one-case}.

\blem\label{le:ordinarity}
Let $E/\Q$ be an elliptic curve and let $p_1<p_2<\cdots<p_m$ be a sequence of primes such that \eqref{eq:prog-seq} holds for some $m\geq 2$. Then in any pair $p_i<p_j$ in this sequence of primes, at most one of $p_i,p_j$ can be a prime of supersingular reduction for $E$. In particular, for any $m\geq 2$, at most one prime $p_1,\ldots,p_m$ in \eqref{eq:prog-seq} can be a prime of  supersingular reduction for $E$.
\elem
\bp 
Suppose that $p\neq q$ are distinct primes of supersingular reduction for $E$ and suppose that $\np{p}=\np{q}$ holds. Then this equality together with supersingularity gives
\be p+1=N_p(E)=N_q(E)=q+1,\ee
which contradicts $p\neq q$. So in any equation $N_p(E)=N_q(E)$ with $p\neq q$, at least one of $p,q$ is a prime of good ordinary reduction for $E$. This proves the assertion.
\ep

\brem 
Let me remark that if $5\leq p<q$ are primes and $\np{p}=\np{q}$, then one has
$$q\leq \frac{\np{p}}{\left(1-\frac{2}{\sqrt{p}}\right)}.$$
In particular, if $ \np{p_1}=\np{p_2}=\cdots=\np{p_m}=N$ provides a chain of primes in an elliptic progression with $5\leq p_1<\cdots< p_m$ then one has $$p_i\leq \frac{N}{\left(1-\frac{2}{\sqrt{p_1}}\right)} \text{ for }i=2,\ldots,m.$$ This remark is useful in searching for primes in elliptic progressions since $N$ is determined by $N=\np{p_1}$. Evidently the first assertion implies the second assertion.
The first assertion is immediate from the Hasse lower bound for $q$:
$$q-2\sqrt{q}\leq q+1-\sqrt{q}\leq \np{q}=\np{p}$$
and as $5\leq p<q$ one has 
$$0<\left(1-\frac{2}{\sqrt{p}}\right)\leq {\left(1-\frac{2}{\sqrt{q}}\right)}$$
and hence
$$q\left(1-\frac{2}{\sqrt{p}}\right)\leq q{\left(1-\frac{2}{\sqrt{q}}\right)}\leq \np{p},$$
which gives the assertion.
\erem

\brem\label{re:CM-k-1-case}
Assume that $E/\Q$ has complex multiplication. In this case by standard results, the CM field $K$ is one of the nine imaginary quadratic fields with class number one and hence the ring of integers is a PID.  Suppose both $p,q$ are primes of ordinary reduction for $E$ (see \Cref{le:ordinarity}). Then both $p,q$ split in $K$. Thus one can write
\be
\begin{aligned} 
N_p(E)&=(\pi-1)(\bar{\pi}-1)\\ 
N_q(E)&=(\tau-1)(\bar{\tau}-1)
\end{aligned}\ee
for some $\pi\in\O_K$ (resp. some $\tau\in\O_K$). Thus we are counting the number of solutions to
$$n=(\pi-1)(\bar{\pi}-1)=(\tau-1)(\bar{\tau}-1)$$
for some elements $\pi,\tau\in K$. This should be considered as \erdos\  type classical problem for $G_2(n)$ in the ring $\O_K$ (\cite{erdos36b}). So I expect that for a CM elliptic curve, the $k=1$ case can be settled using \erdos' method plus some standard number field sieving arguments. Hence I expect that \Cref{con:k-one-case}{\bf(1)} can be settled by methods available to us at the moment. This will be taken up elsewhere.
\erem

\brem 
Comparison of the data for the CM cases (\Cref{49-tab}, \Cref{256-tab}) and the non-CM cases (\Cref{11-tab}, \Cref{43-tab}, \Cref{664-tab}) shows that one should expect different behavior for the counting function $G_1(n)$ in these two cases.
\erem

\brem 
The data in the non-CM cases (\Cref{11-tab}, \Cref{43-tab}, \Cref{664-tab}) shows clear dependence of $G_1(n)$ on Mordell-Weil rank as one might expect because heuristically $$L(E,1)\approx \prod_{p\leq x} \frac{p}{\np{p}}.$$ The dependence on the rank is less clear in the CM cases (\Cref{49-tab}, \Cref{256-tab}), but one should expect it because of this heuristic.
\erem

\section{On the equation $(p_1-1)(p_2-1)\cdots (p_k-1)=n$}
As stated in the \Cref{se:intro}, \cite{erdos36} ($k=3$) and \cite{erdos36b} ($k=2$) were starting points for this paper and I want to remark that the methods I have used in solving \eqref{eq:erdos} for elliptic curves for $k\geq 3$ allows one to solve the general version of Erdos' original problem
\be\label{eq:erdos-orig} 
(p_1-1)(p_2-1)\cdots (p_k-1)=n
\ee
for all $k\geq 2$.
\newcommand{\ekn}{\mathscr{E}_k(n)}
Let $\ekn$ be the number of solutions to \eqref{eq:erdos-orig} in prime numbers $p_1<\cdots < p_k$. 
for all $k\geq2$. The theorem is the following
\bthm\label{th:main2}
For any $k\geq 2$
$$ 
\varlimsup_{n\to\infty} \ekn=\infty.
$$
More precisely, if $k=2$ then for all sufficiently large $x$, there exist  integers $n\leq x$ with $$\ekn\geq e^{\sqrt{\lx}-\vare},$$
and for every integer $k\geq 3$ and for all sufficiently large $x$, there exist  integers $n\leq x$ with $$\ekn\geq (\lx)^{\delta}$$
where  $\delta=\delta(k)=\log(2)k(1-2\vare)-2>0$ (for $k\geq3$).
\ethm

To prove these theorem, one needs the following variants of the results used in the proof of \Cref{th:main2} given earlier.

\bthm\label{th:halberstam}
There exists a positive constant $c_i$ such that for all $x$ sufficiently large one has
$$\sump x (\omega(p-1)-\llx)^2=\cnst x\llx(1+o(1)).$$ 
\ethm
\bp 
This is proved in \cite{halberstam56c}.
\ep

\bpro\label{pro:fourth-moment2}
There is a constant $c_i$ such that 
for all $x$ sufficiently large, one has
$$
\sump x \left(\omega(p-1)-\llx\right)^4=\cnst \pi(x)(\llx)^2+o(\pi(x)(\llx)^2).
$$
\epro
\bp 
This is immediate from \cite[Proposition 3]{granville07}.
\ep

\blem\label{le:basic-count-erdos}
Let $\vare$ be a fixed, sufficiently small positive real number. 
Let $\ax$ denote the number of primes $ p\leq x$ such that 
\be 
\omega(p-1)<\low \ee \
Then there is a positive constant $c_i(\vare)$ such that for all $x\geq x_0$ one has 
$$ 
\cax\leq\cnste \frac{\pi(x)}{(\llx)^2}.
$$
\elem
\bp 
The proof is immediate from the proof of \Cref{le:basic-count1} using \Cref{pro:fourth-moment2} and \Cref{th:halberstam}.
\ep

\blem\label{le:basic-count2-erdos}
Let $\vare$ be a fixed, sufficiently small positive real number. 
Let $\cx$ denote the number of admissible primes i.e. primes $\leq x$ such that 
\be 
\omega(p-1) \geq \low.
\ee \
Then there is a positive constant $c_i(\vare)$ such that for all $x\geq x_0$ one has 
$$ 
\ccx\geq \pi(x)-\cnste\frac{\pi(x)}{(\llx)^2}>\cnst(\vare)\pi(x).
$$
\elem
\bp 
The proof is identical to that of \Cref{le:basic-count2}.
\ep

\bp[Proof of \Cref{th:main2}] 
For $k=2$, I will invoke \cite{erdos36b}. So one has to prove \Cref{th:main2} for $k\geq 3$.
This is carried out using the method of proof of \Cref{th:main} using \Cref{pro:fourth-moment2}, \Cref{le:basic-count-erdos}, \Cref{le:basic-count2-erdos}.

So it will be enough to indicate the quantities which need to be changed in the proof of \Cref{th:main} to arrive at the proof of \Cref{th:main2}. 

Let $\vare$ be a fixed, sufficiently small positive real number. Let $k\geq 1$ be an integer.  Let 
\be 
\bx{k} =\{  p_1\cdot p_2\cdots p_k\leq x:   p_i\neq p_j\ {\rm if}\ i\neq j\  \ {\rm and}\ p_i\ {\rm admissible  \ for\ all}\    1\leq i\leq k\},
\ee
with admissibility in the sense of \Cref{le:basic-count2-erdos}.

Then I claim that it is sufficient to prove that for any $k\geq 1$ one has
\be\label{eq:erdos-bound-orig}
\cbx{k}\geq \cnste\frac{x}{\lx}.
\ee

Suppose for the moment that \eqref{eq:erdos-bound-orig} has been established. Then one can complete the proof of the theorem as follows. On one hand, one has $\cbx{k}\geq \cnste \frac{x}{\lx}$ on the other hand, by \Cref{th:erdos-product} one has $\cdx\leq \cnste\frac{x}{(\lx)^{1+\delta}}$. 

So if one considers the mapping 
\be\label{eq:mapping-orig}
\bx{k}\to \dx\ee 
given by 
\be\label{eq:mapping-orig2} p_1\cdots p_k\in \bx{k} \mapsto (p_1-1)\cdots (p_k-1).\ee 
Then by definition, the number $\ekn$ is precisely the cardinality of the fiber over $n$ under the mapping given by \eqref{eq:mapping-orig}. By \eqref{eq:erdos-bound-orig} and \Cref{th:erdos-product},  some fibers of the mapping \eqref{eq:mapping-orig} must have cardinality $\gg (\log x)^\delta$. Otherwise, for any arbitrary real number $c>0$, the cardinality of any fiber is $\leq c\cdot(\log x )^\delta$. This together with
\Cref{th:erdos-product} (which requires $k\geq3$), gives
\be \cbx{k}\ll \cdx\cdot c\cdot (\log x )^\delta \leq c\cdot \frac{x}{\log(x)},\ee 
which contradicts \eqref{eq:erdos-bound-orig} because $c>0$ is arbitrary. Hence there exist $1\leq n\leq x$ such $G_k(n)\gg (\log x )^\delta$ with $\delta=\delta(k)>0$ given by \eqref{def:deltak}.

So let me now prove the bound asserted in \eqref{eq:erdos-bound-orig}. This is done exactly as in the proof of \Cref{th:main}. For doing this, one needs \Cref{le:series-lemma}. The only point which needs to be changed in the proof of \Cref{le:series-lemma} is that one uses the definition of admissibility given by \Cref{le:basic-count2}. This completes the proof of \Cref{th:main2}.
\ep

\section{Numerical data}\label{se:numerical-data}
\numberwithin{equation}{subsection}
In 2018 when the initial version of this paper (comprising only of Sections \ref{se:intro}--\ref{se:conj}) was   written and submitted for publication, I had very modest evidence for \Cref{con:k-one-case}. After receiving the referee report, and especially because the referee expressed skepticism about \Cref{con:k-one-case}{\bf(1)} (but not about \Cref{con:k-one-case}{\bf(2, 3)}), I decided to gather more evidence for \Cref{con:k-one-case} (I had \Cref{re:CM-k-1-case}, which observes that \Cref{con:k-one-case}{\bf(1)} is amenable by presently available methods in the CM case, at the back of my mind, even though this was not mentioned in my manuscript).   

My initial searches, which led me to formulate \Cref{con:k-one-case}, were limited to a very small data range typically for primes $p\leq 4.36\times 10^8$ (mostly because my modest computer and simple code could not handle higher ranges. The computations get slower for large discriminants.

Eventually, I shifted to a new version of Sage, optimized my code and worked with pre-computed table of primes to extend the search range. This allowed me to expand my computations (gathering raw data  for primes in the range $50\times 10^9 \leq p \leq  70\times 10^9$). This  still took several days run per example on my modest computer.  I have now extensive data $\approx 50{\rm GB}$. However, most of the data remains unanalyzed because during 2020-2025, my work on Mochizuki's work and the $abc$-conjecture came to dominate all other research projects,  and I did not address the revision of this paper requested by the referee until May 2025.

I have now started going over the data and fragments of it are presented below in \Cref{ss:cm-examples} (CM case) and \Cref{ss:non-cm-examples} (non-CM case). The CM behavior is quite distinct from the non-CM case and hence separation of the two cases is essential.  In all cases when I searched in the full range, I found  million(s) of sequences of primes of length at least two i.e. $G_1(n)\geq 2$ (see \Cref{tab:data-for-11a3-3}) and this is irrespective of whether or not the curve has CM. At this point, the evidence for \Cref{con:k-one-case}{\bf(2,3)} is quite compelling.

Long chains of primes are also easy to find for CM-elliptic curves over $\Q$ (\Cref{ss:cm-examples}). Long chains are sparser in the non-CM situation. This is to be expected.  Here are some fragments of the data. \textit{In \Cref{ss:cm-examples}, only the first prime chain of the asserted length is reported.}

\subsection{CM Examples}\label{ss:cm-examples}
The data I have collected clearly separates CM case from the non-CM case. The important difference being the relatively higher probability of finding chains of long (i.e. $>10$) lengths. 
Here are a few examples of this type. 	All of these CM examples were found by a modest brute force search among primes $2\leq p\leq 2036273290\approx 2\times 10^9$. \textit{It should be noted that the list of primes is the first list of primes of the reported length and in the this range.}
\newcommand{\bex}{\begin{example}}
	\newcommand{\eex}{\end{example}}
\bex
For $E:y^2 = x^3 - 27$ (LMFDB label 36.a3)  with CM by $\Q(\zeta_3)$, and $$G_1(327797652)=44.$$ 
Here are the 44 primes with $N_p(E)=327797652$\,:

327761659, 327761719, 327762247, 327762763, 327763507, 327764287, 327766249, 327766657, 327770419, 327771307, 327771739, 327771817, 327776233, 327780139, 327780847, 327784879, 327785083, 327786523, 327788707, 327791017, 327793153, 327795199, 327796939, 327797563, 327802369, 327802987, 327804289, 327806599, 327811753, 327812227, 327815167, 327816343, 327821929, 327823999, 327826003, 327827713, 327828649, 327828967, 327831103, 327831937, 327832513, 327832873, 327833059, 327833833. 
\eex

\bex $E:y^2=x^3-x$ (CM by $\Q(\sqrt{-1})$) (LMFDB label 32.a3)  $$G_1(227232200)=35$$

227202181, 227203829, 227205133, 227207749, 227208181, 227210021, 227211877, 
227213981, 227216413, 227220757, 227222117, 227223221, 227223917, 227226581, 227232037, 227237693, 227237821, 227238533, 227241181, 227242397, 227249837, 227250421, 227252437, 227252621, 227254381, 227254549, 227259581, 227260093, 227260613, 227261189, 227261381, 227261677, 227262197, 227262293, 227262349.
\eex 

\bex The elliptic curve $E:y^2 + y = x^3$ (LMFDB label 27.a4). This has CM by $\Q(\sqrt{-3})$. Here are the $43$ primes $p$ for which $$G_1(1530144252)=43$$

1530066103, 1530068047, 1530069763, 1530070933, 1530073327, 1530074983,  1530075847, 1530080059, 1530083839, 1530083959 1530090823, 1530102037, 
1530105019, 1530105253, 1530106603, 1530113989, 1530118159, 1530121477, 
1530127717, 1530132943, 1530142927, 1530145579, 1530147127, 1530147397, 
1530151123, 1530165283, 1530167029, 1530170347, 1530174739, 1530179557, 
1530180187, 1530181903, 1530191407, 1530191599, 1530191749, 1530202207, 
1530208447, 1530210379, 1530210523, 1530216949, 1530221299, 1530221893,  1530222469.
\eex
\bex
\item For $y^2=x^3-35x+98$ (LMFDB label 784.f4) with CM by $\Q(\sqrt{-7})$, here is a chain with $$G_1(273758144)=28$$

273725797, 273726443, 273726721, 273727871, 273734003, 273737333, 273739069, 273744277,
273746519, 273751783, 273753547, 273755393, 273756737, 273757663, 273766819, 273768133,
273771521, 273773677, 273775919, 273777793, 273779669, 273782021, 273782107, 273783617,
273785717, 273786221, 273789077, 273791207. 
\eex

\begin{example}
In {\sc\Cref{256-tab}} one finds a more detailed table for the CM elliptic curve $ y^2 = x^{3} - 2 x $ (LMFDB label 256.b1) in the small search range $ 257 \leq p \leq 430000000 $, the number of of  primes in this range = 19617094. As one sees from {\sc\Cref{256-tabB}}, primes of  longer length occur quite frequently in the CM case and  prime chains of length $\geq 20$ are listed in {\sc\Cref{256-tabA}}.
\end{example}

{\small 
	\afterpage{
\newgeometry{margin=0.5in}
\thispagestyle{empty}
\begin{landscape}
 {\tiny \begin{table}[!htb] 
 \caption{{\tiny Elliptic curve $ y^2 + x y = x^{3} -  x^{2} - 2 x - 1 $, CM, rank zero, LMFDB label 49.a4 
 search range = [ 50 $\leq p \leq$ 430000000 ] primes in range = 18083595 } }\label{49-tab} 
 \begin{subtable}[t]{.75\linewidth}
\flushleft\caption{ prime chains }
\begin{tabular}{|p{0.1125\linewidth}|p{0.67\linewidth}|} \hline 
  n & prime chain (lengths $\geq 29$) \\ 
 \hline132159104 &  132136117, 132136159, 132136313, 132136943, 132138553, 132140009, 132140933, 132141577, 132144013, 132144559, 132148787, 132151993, 132153253, 132156067, 132158951, 132159103, 132162143, 132164957, 132166217, 132168233, 132171901, 132173651, 132174197, 132176633, 132178901, 132179657, 132179909, 132181267, 132181337  \\ \hline 
270389504 &  270356633, 270356717, 270356857, 270357907, 270359909, 270360553, 270361001, 270366601, 270367637, 270369667, 270370409, 270371837, 270373601, 270376793, 270381497, 270386117, 270386383, 270387937, 270389351, 270389659, 270392893, 270396617, 270402217, 270404359, 270408601, 270412409, 270418009, 270418457, 270421103, 270421733, 270421901, 270422153  \\ \hline 
306403328 &  306368497, 306368609, 306369169, 306369617, 306370457, 306370807, 306371801, 306374447, 306379529, 306380257, 306384289, 306387271, 306387509, 306391457, 306402461, 306406927, 306416909, 306419149, 306419387, 306420017, 306426569, 306427129, 306430097, 306434017, 306434521, 306434857, 306435851, 306435949, 306436201  \\ \hline 
\end{tabular} 
 \end{subtable}

\begin{subtable}[t]{.25\linewidth}
\flushright
\caption{Chain counts per length are }
 \begin{tabular}{|c|c|} \hline length & number of prime chains \\ \hline 
2 & 1569802 \\ \hline 
3 & 593954 \\ \hline 
4 & 248637 \\ \hline 
5 & 114945 \\ \hline 
6 & 56577 \\ \hline 
7 & 29581 \\ \hline 
8 & 16422 \\ \hline 
9 & 9233 \\ \hline 
10 & 5376 \\ \hline 
11 & 3220 \\ \hline 
12 & 2006 \\ \hline 
13 & 1252 \\ \hline 
14 & 803 \\ \hline 
15 & 532 \\ \hline 
16 & 320 \\ \hline 
17 & 197 \\ \hline 
18 & 154 \\ \hline 
19 & 97 \\ \hline 
20 & 67 \\ \hline 
21 & 46 \\ \hline 
22 & 17 \\ \hline 
23 & 14 \\ \hline 
24 & 15 \\ \hline 
25 & 8 \\ \hline 
26 & 7 \\ \hline 
27 & 8 \\ \hline 
28 & 4 \\ \hline 
29 & 2 \\ \hline 
32 & 1 \\ \hline 
\end{tabular}
 \end{subtable}
 \end{table} 
 }
\end{landscape}
\clearpage
\restoregeometry
}
 }
\newpage
{\small 
	\afterpage{
\newgeometry{margin=0.5in}
\thispagestyle{empty}
\begin{landscape}
 {\tiny \begin{table}[!htb] 
 \caption{{\tiny CM Elliptic curve $ y^2 = x^{3} - 2 x $ CM, rank one, LMFDB label 256.b1 
 search range = [ 257 $\leq p \leq$ 430000000 ] primes in range = 19617094 } }\label{256-tab} 
 \begin{subtable}{\linewidth}
\centering\caption{\tiny prime chains of lengths $\geq 20$  }\label{256-tabA}
\begin{tabular}{|p{0.07\linewidth}|p{0.9\linewidth}|} \hline 
  n & prime chain \\ 
 \hline 130832000 & {\tiny 130809137, 130809233, 130809361, 130809841, 130810129, 130810481, 130811921, 130815121, 130830577, 130831249, 130832753, 130841489, 130842097, 130842961, 130846321, 130852081, 130853137, 130853521, 130854161, 130854481, 130854769, 130854833} \\ \hline 
178568000 & {\tiny 178541281, 178541329, 178543121, 178544033, 178544161, 178549729, 178553953, 178555921, 178560881, 178561057, 178567999, 178569697, 178582049, 178585361, 178588321, 178589041, 178591969, 178592369, 178593809, 178594321, 178594721 } \\ \hline 
183010880 & {\tiny 182983873, 182985329, 182988433, 182991377, 182992129, 182995777, 182998513, 183006337, 183006529, 183010753, 183010879, 183011009, 183015233, 183021169, 183025793, 183028337, 183031553, 183033457, 183036097, 183037073, 183037889 } \\ \hline 
217393280 & {\tiny 217363793, 217363889, 217364113, 217365233, 217366129, 217368337, 217370353, 217372049, 217379249, 217381553, 217384177, 217386193, 217393457, 217396913, 217405009, 217407313, 217409009, 217410257, 217416209, 217417393, 217418897, 217420433 } \\ \hline 
258128000 & {\tiny 258096017, 258100561, 258101297, 258101873, 258103441, 258105041, 258107281, 258109297, 258117361, 258123793, 258127999, 258131089, 258132209, 258132881, 258143153, 258145777, 258148721, 258150481, 258152561, 258154769, 258155441, 258157841, 258158321, 258159761, 258159857 } \\ \hline 
266614400 & {\tiny 266581969, 266582161, 266582801, 266585489, 266586161, 266591729, 266594641, 266602033, 266605457, 266606161, 266614399, 266614609, 266622001, 266622641, 266630801, 266631601, 266640401, 266644817, 266646833, 266646929 } \\ \hline 
276360500 & {\tiny 276327361, 276327593, 276328681, 276328961, 276332033, 276338473, 276347881, 276355841, 276357161, 276359233, 276362921, 276371009, 276372353, 276373121, 276386497, 276387841, 276388969, 276391657, 276392041, 276392321, 276393281, 276393409 } \\ \hline 
279012500 & {\tiny 278979097, 278979161, 278979601, 278980241, 278981401, 278984081, 278988121, 278990801, 278992081, 278997721, 279012499, 279013201, 279021401, 279027601, 279036881, 279039641, 279042961, 279043601, 279045401, 279045841 } \\ \hline 
327080000 & {\tiny 327044033, 327044401, 327045041, 327046081, 327046961, 327049601, 327051601, 327060401, 327067441, 327069697, 327076177, 327077441, 327092561, 327093553, 327097601, 327102401, 327104641, 327105409, 327112321, 327113921, 327115601, 327116161 } \\ \hline 
348296000 & {\tiny 348258721, 348259361, 348263009, 348267281, 348268049, 348268321, 348272161, 348274033, 348275489, 348280081, 348304673, 348306961, 348308641, 348313457, 348316513, 348317969, 348320737, 348323681, 348324721, 348328993, 348331681, 348332641 } \\ \hline 
420430400 & {\tiny 420391441, 420391553, 420394097, 420396433, 420398833, 420399361, 420406001, 420417601, 420418673, 420421361, 420421649, 420430657, 420434689, 420439153, 420449009, 420451441, 420453377, 420456577, 420456977, 420461441, 420464017, 420471409 } \\ \hline 
425444500 & {\tiny 425403353, 425406361, 425407769, 425409041, 425410897, 425412121, 425416729, 425418457, 425419409, 425434193, 425436761, 425438681, 425458841, 425473681, 425479513, 425482961, 425484049, 425484377, 425485337, 425485649 } \\ \hline 
\end{tabular} 
 \end{subtable}
\begin{subtable}{.5\linewidth}
 \centering
 \caption{\tiny Chain counts per length are }\label{256-tabB}
 \begin{tabular}{|c|c|} \hline length & number of prime chains \\ \hline 
2 & 1584489 \\ \hline 
3 & 407677 \\ \hline 
4 & 130577 \\ \hline 
5 & 49690 \\ \hline 
6 & 20715 \\ \hline 
7 & 9525 \\ \hline 
8 & 4620 \\ \hline 
9 & 2351 \\ \hline 
10 & 1284 \\ \hline 
11 & 678 \\ \hline 
12 & 377 \\ \hline 
13 & 190 \\ \hline 
14 & 112 \\ \hline 
15 & 67 \\ \hline 
16 & 35 \\ \hline 
17 & 23 \\ \hline 
18 & 16 \\ \hline 
19 & 7 \\ \hline 
20 & 3 \\ \hline 
21 & 2 \\ \hline 
22 & 6 \\ \hline 
25 & 1 \\ \hline 
\end{tabular}
 \end{subtable}
 \end{table} 
 }
\end{landscape}
\clearpage
\restoregeometry
}

 }
\clearpage
\newpage

\subsection{Non-CM examples}\label{ss:non-cm-examples}
\begin{example}  Let $E/\Q$ be the elliptic curve 
	$E:y^2 + 3 y = x^{3} -  x + 2$ (LMFDB label 7739.a1). The {\sc\Cref{7739-tabA}} gives a small list $E$-prime triples i.e. solutions to the equation
	$G_1(n)=3$ (in other words $n$ for which there are there primes $p_1,p_2,p_3$ with $\np{p_i}=n$);
	  {\sc\Cref{7739-tabB}} provides a number of prime chains found in the range stated at the top of the table; {\sc\Cref{7739-tabC}} provides a small table of values of the function $G_1(n)$.
	  
\end{example}

\brem 
Examining  {\sc\Cref{7739-tabA}}, it is clear that the triples of primes listed about need not form an arithmetic progression of primes. For example for prime triples in the $2^{nd}, 3^{rd}$ and $4^{th}$ rows, corresponding to $n=11422,12312,12672$, the successive differences are
\be
\begin{aligned}
11519-11299 &= 2^2\times 5 \times 11\\
11617-11519 &=2\times 7^2\\
12391-12161 & = 2\times 5\times 23\\
12421-12391 & = 2\times 5\times 3\\
12721-12619 & = 2\times 3\times 17\\
12791-12721 & = 2\times 5 \times 7.
\end{aligned}
\ee	
These equations show the each of the first three $E$-triplets of primes in the above table  do not form an arithmetic progression of primes. Hence one is dealing with the occurrence of new phenomenon in primes in the presence of the elliptic curve $E/\Q$.
\erem
 {\tiny \begin{table}[!htb] 
 \caption{{\tiny Elliptic curve $ y^2 + y = x^{3} -  x + 4 $ LMFDB label 7739.a1 
 search range = [ 7740 $\leq p \leq$ 50000 ] primes in range = 3906 } } 
 \begin{subtable}{0.33\linewidth}
\centering\subcaption{ prime chains }\label{7739-tabA}
\begin{tabular}{|c|l|} \hline 
  n & prime chain \\ 
 \hline8568 & $ 8423 , 8527 , 8563 $ \\ \hline 
11422 & $ 11299 , 11519 , 11617 $ \\ \hline 
12312 & $ 12161 , 12391 , 12421 $ \\ \hline 
12672 & $ 12619 , 12721 , 12791 $ \\ \hline 
32022 & $ 31699 , 31873 , 32213 $ \\ \hline 
34240 & $ 34217 , 34327 , 34603 $ \\ \hline 
37464 & $ 37517 , 37571 , 37693 $ \\ \hline 
46992 & $ 46691 , 46747 , 46817 $ \\ \hline 
48528 & $ 48353 , 48487 , 48563 $ \\ \hline 
\end{tabular} 
 \end{subtable}
\begin{subtable}{0.33\linewidth}
 \centering
 \caption{Chain counts per length are }\label{7739-tabB}
 \begin{tabular}{|c|c|} \hline length & number of prime chains \\ \hline 
2 & 228 \\ \hline 
3 & 9 \\ \hline 
\end{tabular}
 \end{subtable}
 \begin{subtable}{0.33\linewidth}
 \caption{{\tiny Small table of values of $G_1(n)\geq 2$}}\label{7739-tabC}
 	\begin{tabular}{|c|c|}
 \hline $n$ & $G_1(n)$ \\ 
 \hline 
 \hline 10262 & 2\\
 \hline 10494 & 2\\
 \hline 10630 & 2\\
 \hline 10697 & 2\\
 \hline 10704 & 2\\
 \hline 11072 & 2\\
 \hline 11100 & 2\\
 \hline 11168 & 2\\
 \hline 11276 & 2\\
 \hline 11422 & 3\\
 \hline 11441 & 2\\
 \hline
 \end{tabular} 
 \end{subtable}
 \end{table} 
 } \newpage

\begin{example}[The modular curve $X_1(11)$ (LMFDB label 11.a3)]
Let $E: y^2+y=x^3-x^2$. This is the genus one modular curve $X_1(11)$, with Mordell-Weil rank zero, and with  LMFDB label 11.a3. The data I have collected for this is extensive {\sc\Cref{11-tab}} gives a list prime chains in a relative small range $12\leq p\leq 43\times 10^6$. 

{\sc\Cref{tab:data-for-11a3}}, {\sc\Cref{tab:data-for-11a3-2}}, and {\sc\Cref{tab:data-for-11a3-3}}  give data (for this curve) for primes in a substantially larger block $50\times 10^9 \leq p\leq 65\times 10^9$. {\sc\Cref{tab:data-for-11a3}}, {\sc\Cref{tab:data-for-11a3-2}} list some prime chains,  while {\sc\Cref{tab:data-for-11a3-3}} lists data on the number of chains in this range.
	
Here are some highlights of {\sc\Cref{11-tab}} and {\sc\Cref{tab:data-for-11a3-3}}. These tables refer to sets of search ranges differing by an order of magnitude. In both the ranges the longest chain is of length 10. In the higher range of {\sc\Cref{tab:data-for-11a3-3}}, the chain of length ten corresponds to	 $n=57346191000$ i.e. one has $G_1(57346191000)=10$. The ten primes $p$ with $$N_p(E)=57346191000$$
are listed below:
\\ \\
{\small
	\indent 57345977191, 57346023283, 57346109347, 57346128311, 57346137419, 57346254151, 57346262929, 57346303661, 57346391807, 57346523111. All these primes are $\approx 57\times 10^9$.}
\\ 
	\\
	A purely randomized, and highly non-deterministic and computationally expensive, search  finds chains in a substantially larger range for this curve. One has
	$$G_1(7816253637600)=2$$ and here are the two primes ($\approx 7.8\times 10^{12}$) for which this happens $$7816251016289, 7816252108193.$$
	Note that $$7816252108193-7816251016289=1091904.$$
{\small 
	\begin{landscape}
\begin{center}
 {\tiny \begin{table}[!htb] 
 \caption{{\tiny Elliptic curve $ y^2 + y = x^{3} -  x^{2} $, rank zero, LMFDB label 11.a3 
 search range = [ 12 $\leq p \leq$ 430000000 ] primes in range = 18920234 } }\label{11-tab} 
 \begin{subtable}{\linewidth}
\centering\caption{ prime chains (lengths $\geq 9$)}\label{11-tabA}
\begin{tabular}{|c|l|} \hline 
  n & prime chain \\ 
 \hline508200 & $ 507151 , 507631 , 507827 , 507971 , 508021 , 508271 , 508301 , 508619 , 508817 $ \\ \hline 
5704000 & $ 5701651 , 5702003 , 5702581 , 5703001 , 5703521 , 5704081 , 5704201 , 5707459 , 5707721 $ \\ \hline 
19501200 & $ 19495031 , 19495291 , 19495981 , 19496941 , 19499671 , 19501501 , 19502311 , 19504913 , 19506367 $ \\ \hline 
20946600 & $ 20940131 , 20942741 , 20944087 , 20944291 , 20944571 , 20947807 , 20948003 , 20951291 , 20951971 $ \\ \hline 
28987200 & $ 28980191 , 28982101 , 28983391 , 28984009 , 28984751 , 28985549 , 28989211 , 28991981 , 28992991 $ \\ \hline 
31693200 & $ 31683721 , 31685461 , 31692191 , 31692851 , 31697929 , 31700239 , 31701701 , 31701941 , 31702087 $ \\ \hline 
34272000 & $ 34260761 , 34267991 , 34268621 , 34269269 , 34275211 , 34275671 , 34276813 , 34279961 , 34281461 $ \\ \hline 
53982600 & $ 53970487 , 53975767 , 53978189 , 53980481 , 53980981 , 53985431 , 53985823 , 53991719 , 53994203 $ \\ \hline 
86293200 & $ 86276051 , 86284931 , 86285891 , 86286817 , 86291761 , 86297741 , 86297903 , 86302039 , 86307343 $ \\ \hline 
147237600 & $ 147225601 , 147228377 , 147238601 , 147238781 , 147238981 , 147241753 , 147242801 , 147249281 , 147257483 $ \\ \hline 
204745200 & $ 204717221 , 204720097 , 204727093 , 204735031 , 204739441 , 204743711 , 204744343 , 204749221 , 204756401 , 204758137 $ \\ \hline 
232312500 & $ 232287089 , 232287859 , 232293821 , 232294421 , 232309201 , 232318921 , 232319641 , 232324531 , 232329107 $ \\ \hline 
235407600 & $ 235400153 , 235404061 , 235404671 , 235409011 , 235416127 , 235420291 , 235422541 , 235423849 , 235424617 $ \\ \hline 
286923600 & $ 286900771 , 286911967 , 286917391 , 286921889 , 286922611 , 286924061 , 286926509 , 286936913 , 286941269 $ \\ \hline 
\end{tabular} 
 \end{subtable}
\begin{subtable}{.5\linewidth}
 \centering
 \vskip0.5in
 \caption{Chain counts per length are }
 \begin{tabular}{|c|c|} \hline length & number of prime chains \\ \hline 
2 & 2686888 \\ \hline \label{11-tabB}
3 & 458962 \\ \hline 
4 & 82371 \\ \hline 
5 & 14768 \\ \hline 
6 & 2690 \\ \hline 
7 & 466 \\ \hline 
8 & 65 \\ \hline 
9 & 13 \\ \hline 
10 & 1 \\ \hline 
\end{tabular}
 \end{subtable}
 \end{table} 
 }
\end{center}
 \end{landscape} }
\newpage	
\thispagestyle{empty}
\clearpage
	\begin{landscape}
\begin{table}
		\caption{Fragment of data for $E:y^2+y=x^3-x^2$ i.e. the modular curve $X_{1}(11)$ (LMFDB label 11.a3)}\label{tab:data-for-11a3}
\begingroup\scriptsize	
		\setlength{\tabcolsep}{10pt} \renewcommand{\arraystretch}{2}
		\begin{tabular}{|c|c|c|c|c|c|c|c|c|c|}
			\hline 
			\rule[-1ex]{0pt}{2.5ex} $n$ & \multicolumn{8}{c|}{ complete list of primes $p$ such that $N_p(E)=n$ (the primes $p$ are $\approx 57\times 10^9$)} \\
			\hline
			\rule[-1ex]{0pt}{2.5ex}
			57273393600 &57273249581  & 57273266863 & 57273280231 & 57273451451 & 57273474851 & 57273595331 & 57273619201 & \\
			\hline
			\rule[-1ex]{0pt}{2.5ex} 57256718800 &57256478821 & 57256508581 & 57256667441 & 57256727659 & 57256799011 & 57256849387 & 57257141341 & \\
			\hline
			\rule[-1ex]{0pt}{2.5ex} 57257048800 &57256662901 & 57256945307 & 57257072911 & 57257256167 & 57257421409 & 57257458691 & 57257464051 & 57257484281\\
			\hline
			\rule[-1ex]{0pt}{2.5ex} 57286867500 & 57286629791 & 57286697071 & 57286762901 & 57286768087 & 57286796753 & 57286922681 & 57287137601 & \\
			\hline
			\rule[-1ex]{0pt}{2.5ex} 57279320400 &57278958031 & 57278995819 & 57279251851 & 57279276713 & 57279528569 & 57279679361 & 57279691019 & \\
			\hline
			\rule[-1ex]{0pt}{2.5ex} 57266077440 &57265778333 & 57265940143 & 57266028377 & 57266038033 & 57266278847 & 57266371187 & 57266393297 & \\
			\hline
			\rule[-1ex]{0pt}{2.5ex} 57271311000 &57270891101 & 57271097291 & 57271161061 & 57271173961 & 57271251491 & 57271361179 & 57271456901 & \\
			\hline
			\rule[-1ex]{0pt}{2.5ex} 57275666400 &57275247031 & 57275488117 & 57275608171 & 57275618819 & 57275629049 & 57275679053 & 57276099821 & \\
			\hline
			\rule[-1ex]{0pt}{2.5ex} 57280327000 &57279891901 & 57279995057 & 57280141861 & 57280223551 & 57280244317 & 57280591831 & 57280600891 & \\
			\hline
			\rule[-1ex]{0pt}{2.5ex} 57284579400 &57284432911 & 57284516761 & 57284620451 & 57284657341 & 57284767297 & 57285027871 & 57285035477 & \\
			\hline
			\rule[-1ex]{0pt}{2.5ex} 57277954800 & 57277700743 & 57277733041 & 57277989851 & 57278008259 & 57278056061 & 57278127983 & 57278278639 & \\
			\hline
			\rule[-1ex]{0pt}{2.5ex} 57268260000 &57268061311 & 57268082581 & 57268198307 & 57268465627 & 57268471451 & 57268472293 & 57268715767 & \\
			\hline
			\rule[-1ex]{0pt}{2.5ex}57293937600 &57293498611 & 57293579611 & 57293671811 & 57293714371 & 57294096397 & 57294181291 & 57294233897 & \\
			\hline
		\end{tabular}
\endgroup 
		\end{table}	
	\end{landscape}

	\clearpage	
	\newpage
	\thispagestyle{empty}

	\begin{landscape}
		\begingroup\scriptsize	
		\setlength{\tabcolsep}{10pt} \renewcommand{\arraystretch}{2}
		\begin{table}
			\caption{Searching prime chains for primes between 50000000000 and 65000000000 for $E:y^2+y=x^3-x^2$ i.e. the modular curve $X_{1}(11)$ (LMFDB label 11.a3)}
			\label{tab:data-for-11a3-2}
			\begin{tabular}{|c|c|l|}
				\hline 
				chain length & $n$ & all the primes with $50\times 10^9 \leq p \leq 65\times 10^9$ and  $N_p(E)=n$ (table shows the first chain of the reported length) \\
				\hline 
				2 & 50166672045 & 50166834877 50166962849 \\
				\hline 
				3 & 50418330320 & 50418217649, 50418221993, 50418515317\\
				\hline
				4 & 50314871100 & 50314564921, 50314618051, 50314946659, 50315158801\\
				\hline
				5 & 50130321600 & 50129912909, 50130144851, 50130201271, 50130706999, 50130731801\\ 
				\hline
				6 & 50482649100 & 50482313573, 50482454923, 50482493557, 50482729439, 50482879327, 50482881989\\
				\hline
				7 & 50499667800 & 50499444929, 50499493651, 50499649861, 50499785291, 50499803791, 50499817561, 50499872701\\ 
				\hline
				8 & 50433188400 & 50432916953, 50432917511, 50433029173, 50433105757, 50433310181, 50433411181, 50433465491, 50433502291\\ 
				\hline
				9 & 50371653600 & 50371277263, 50371291121, 50371552991, 50371576331, 50371682129, 50371713001, 50371731281, 50371742893, 50371812391\\ 
				\hline
				10 & 57346191000 & 57345977191, 57346023283, 57346109347, 57346128311, 57346137419, 57346254151, 57346262929, 57346303661, 57346391807, 57346523111\\ 
				\hline
			\end{tabular}
		\end{table}
		\vskip2mm 
		\endgroup 	
\end{landscape}
	
	\clearpage	
	\newpage
	\thispagestyle{empty}
	\begin{landscape}
		\begingroup\scriptsize	
		\setlength{\tabcolsep}{10pt} \renewcommand{\arraystretch}{2}
		\begin{table}
			\caption{Number of prime chains with $50\times 10^9 \leq p\leq 65\times 10^9$ for $E:y^2+y=x^3-x^2$ i.e. the modular curve $X_{1}(11)$ (LMFDB label 11.a3)}
			\label{tab:data-for-11a3-3}
			\begin{tabular}{|c|c|}
				\hline
				length & Number of chains of each length with $50\times 10^9 \leq p\leq 65\times 10^9$ \\
				\hline
				2 & 38965829\\
				\hline
				3 & 5269795\\
				\hline
				4 & 753485\\
				\hline
				5 & 107577\\
				\hline
				6 & 14901\\
				\hline
				7 & 1829\\
				\hline
				8 & 253\\
				\hline
				9 & 21\\
				\hline
				10 & 1\\
				\hline
			\end{tabular}
		\end{table}
		\vskip2mm 
		\endgroup 	
\end{landscape}
\end{example}

{\small 
	 {\tiny \begin{table}[!htb] 
 \caption{{\tiny Elliptic curve $ y^2 + y = x^{3} + x^{2} $, rank one, LMFDB label 43.a1 
 search range = [ 44 $\leq p \leq$ 430000000 ] primes in range = 22099853 } }\label{43-tab} 
 \begin{subtable}{.5\linewidth}
\centering\caption{the $16$ prime chains of length $5$}
\begin{tabular}{|c|l|} \hline 
  n & prime chain \\ 
 \hline2275810 & $ 2274487 , 2275219 , 2275439 , 2275457 , 2276137 $ \\ \hline 
33380440 & $ 33377779 , 33379943 , 33381629 , 33383069 , 33384677 $ \\ \hline 
50984136 & $ 50985481 , 50986589 , 50992153 , 50992441 , 50996879 $ \\ \hline 
68962968 & $ 68959043 , 68960503 , 68962379 , 68969501 , 68974063 $ \\ \hline 
137003976 & $ 136986709 , 136992533 , 136993931 , 137000837 , 137019937 $ \\ \hline 
150592440 & $ 150568657 , 150575189 , 150592583 , 150592781 , 150609289 $ \\ \hline 
154869900 & $ 154851233 , 154851757 , 154867463 , 154873547 , 154881277 $ \\ \hline 
171995824 & $ 171972499 , 171981851 , 172000949 , 172018591 , 172019693 $ \\ \hline 
182774124 & $ 182755093 , 182768191 , 182772671 , 182793491 , 182799821 $ \\ \hline 
206888904 & $ 206870039 , 206874959 , 206882827 , 206883407 , 206912179 $ \\ \hline 
251860752 & $ 251845229 , 251848501 , 251869477 , 251875313 , 251886923 $ \\ \hline 
345092372 & $ 345068527 , 345070963 , 345077591 , 345089881 , 345094721 $ \\ \hline 
346218360 & $ 346190003 , 346191029 , 346201421 , 346228889 , 346237139 $ \\ \hline 
350316828 & $ 350284453 , 350286103 , 350290093 , 350341631 , 350344937 $ \\ \hline 
365565420 & $ 365550149 , 365551649 , 365572313 , 365583677 , 365596177 $ \\ \hline 
388397520 & $ 388390699 , 388395677 , 388405907 , 388425881 , 388430629 $ \\ \hline 
\end{tabular} 
 \end{subtable}
\begin{subtable}{.5\linewidth}
 \centering
 \caption{Chain counts per length are }
 \begin{tabular}{|c|c|} \hline length & number of prime chains \\ \hline 
2 & 707341 \\ \hline 
3 & 19624 \\ \hline 
4 & 510 \\ \hline 
5 & 16 \\ \hline 
\end{tabular}
 \end{subtable}
 \end{table} 
 } }
\newpage

{\small 
	 {\tiny \begin{table}[!htb] 
 \caption{{\tiny Elliptic curve $ y^2 = x^{3} - 7 x + 10 $, rank two, LMFDB label 664.a1 
 search range = [ 665 $\leq p \leq$ 430000000 ] primes in range = 22101208 } }\label{664-tab}
 \begin{subtable}{.5\linewidth}
\centering\caption{ prime chains of length $\geq 5$ }
\begin{tabular}{|c|l|} \hline 
  n & prime chain \\ 
 \hline7197528 & $ 7193093 , 7193509 , 7196369 , 7199057 , 7199639 $ \\ \hline 
11064240 & $ 11059339 , 11060173 , 11066089 , 11068033 , 11068289 $ \\ \hline 
36745200 & $ 36741469 , 36741919 , 36745133 , 36746653 , 36754079 , 36754217 $ \\ \hline 
94956516 & $ 94944149 , 94948459 , 94949969 , 94950929 , 94958779 $ \\ \hline 
177483934 & $ 177474931 , 177480367 , 177490171 , 177503503 , 177506711 $ \\ \hline 
183767952 & $ 183748049 , 183763451 , 183774211 , 183785333 , 183786179 $ \\ \hline 
186600656 & $ 186592487 , 186598327 , 186605513 , 186610057 , 186620843 $ \\ \hline 
194959380 & $ 194948059 , 194967007 , 194971033 , 194973959 , 194984039 $ \\ \hline 
296511280 & $ 296484299 , 296499253 , 296528663 , 296531531 , 296538029 $ \\ \hline 
319415616 & $ 319383821 , 319407971 , 319415861 , 319423283 , 319436813 $ \\ \hline 
348706836 & $ 348679717 , 348689993 , 348708229 , 348728669 , 348738707 $ \\ \hline 
401281444 & $ 401262809 , 401263151 , 401285779 , 401287741 , 401307737 $ \\ \hline 
\end{tabular} 
 \end{subtable}
\begin{subtable}{.5\linewidth}
 \centering
 \caption{Chain counts per length are }
 \begin{tabular}{|c|c|} \hline length & number of prime chains \\ \hline 
2 & 706558 \\ \hline 
3 & 19328 \\ \hline 
4 & 486 \\ \hline 
5 & 11 \\ \hline 
6 & 1 \\ \hline 
\end{tabular}
 \end{subtable}
 \end{table} 
 } }
\newpage

\bibliographystyle{plainnat}

 \end{document}